\documentclass[11pt,reqno,a4paper]{article} 

\usepackage{anysize}
\marginsize{3.0cm}{3.0cm}{2.2cm}{2.2cm}

\usepackage[english]{babel}
\usepackage[centertags]{amsmath}
\usepackage{amsfonts,amssymb,amsthm}
\usepackage{color}
\usepackage{hyperref}
\usepackage{mathrsfs}
\usepackage[sans]{dsfont}
\usepackage{tcolorbox}

\let\OLDthebibliography\thebibliography
\renewcommand\thebibliography[1]{
  \OLDthebibliography{#1}
  \setlength{\parskip}{0pt}
  \setlength{\itemsep}{0pt plus 0.3ex} }

\numberwithin{equation}{section}

\theoremstyle{plain}
\newtheorem{theorem}{Theorem}[section]
\newtheorem{proposition}[theorem]{Proposition}

\theoremstyle{definition}
\newtheorem{remarkbase}[theorem]{Remark}
\newenvironment{remark}{\pushQED{\qed} \remarkbase}{\popQED\endremarkbase}

\newcommand{\N}{{\mathbb N}}
\newcommand{\R}{{\mathbb R}}
\newcommand{\C}{{\mathbb C}}
\newcommand{\Z}{\mathbb Z}

\newcommand{\T}{{\mathbb T}}

\newcommand{\mA}{\mathcal{A}}
\newcommand{\mB}{\mathcal{B}}
\newcommand{\mC}{\mathcal{C}}

\newcommand{\mE}{\mathcal{E}}
\newcommand{\mF}{\mathcal{F}}
\newcommand{\mG}{\mathcal{G}}
\newcommand{\mH}{\mathcal{H}}

\newcommand{\mL}{\mathcal{L}}

\newcommand{\mN}{\mathcal{N}}

\renewcommand{\a}{\alpha}
\renewcommand{\b}{\beta}
\newcommand{\g}{\gamma}
\renewcommand{\d}{\delta}

\newcommand{\e}{\varepsilon}
\newcommand{\ph}{\varphi}
\newcommand{\lm}{\lambda}

\newcommand{\om}{\omega}

\newcommand{\s}{\sigma}
\renewcommand{\t}{\tau}

\newcommand{\ttf}{\mN}

\newcommand{\gr}{\nabla}

\newcommand{\pa}{\partial}

\newcommand{\bcr}{\begin{color}{red}}
\newcommand{\ec}{\end{color}}
\newcommand{\bcb}{\begin{color}{blue}}
\newcommand{\bcg}{\begin{color}{green}}

\renewcommand{\oe}{\"o}
\newcommand{\Hor}{H\"ormander}

\title{A Nash-Moser-H\"ormander implicit function theorem\\ 
with applications to control and Cauchy problems for PDEs}

\author{\small{Pietro Baldi, Emanuele Haus}}
\date{} 
\begin{document}

\maketitle

\begin{small}
\textbf{Abstract.}
We prove an abstract Nash-Moser implicit function theorem 
which, when applied to control and Cauchy problems for PDEs in Sobolev class, 
is sharp in terms of the loss of regularity of the solution of the problem 
with respect to the data. 
The proof is a combination of: $(i)$ the iteration scheme by H\oe{}rmander (ARMA 1976),
based on telescoping series, and very close to the original one by Nash;  
$(ii)$ a suitable way of splitting series in scales of Banach spaces, 
inspired by a simple, clever trick used in paradifferential calculus 
(for example, by M\'etivier).  
As an example of application, 
we apply our theorem to a control and a Cauchy problem 
for quasi-linear perturbations of KdV equations, 
improving the regularity of a previous result.
With respect to other approaches to control and Cauchy problems, 
the application of our theorem requires lighter assumptions to be verified.
\emph{MSC2010:} 47J07, 35Q53, 35Q93.
\end{small}

\medskip

\emph{Contents.}
\ref{sec:intro} Introduction ---
\ref{sec:NM} {A Nash-Moser-\Hor{} theorem} ---
\ref{sec:proof} Proof of Theorem \ref{thm:1} ---
\ref{sec:appl} Application to quasi-linear perturbations of KdV.

\section{Introduction}
\label{sec:intro}

In this paper we prove an abstract Nash-Moser implicit function theorem 
(Theorem \ref{thm:1}) which, 
when applied to control and Cauchy problems for evolution PDEs in Sobolev class, 
is sharp in terms of the loss of regularity of the solution of the problem 
with respect to the data.

In terms of such a loss, the sharpest Nash-Moser theorem in literature 
seems to be the one by H\oe{}rmander 
(Theorem 2.2.2 in Section 2.2 of \cite{Geodesy}, and main Theorem in \cite{Olli}). 
\Hor's theorem is sharp when applied to PDEs in H\oe lder spaces (with non-integer exponent), 
but it is \emph{almost} sharp in Sobolev class: if the approximate right inverse 
of the linearized operator loses $\g$ derivatives, and the data of the problem belong to $H^s$, 
then the application of \Hor's theorem gives solutions of regularity $H^{s - \g - \e}$ 
for all $\e > 0$, whereas one expects to find $H^{s - \g}$ 
(and in many cases, with other techniques, 
in fact one can prove such a sharp regularity).
Our Theorem \ref{thm:1} applies to Sobolev spaces with sharp loss,
and thus it extends \Hor's result to Sobolev spaces.

\medskip

As it is well-known, the Nash-Moser approach is 
natural to use in situations where a loss of regularity 
prevents the application of other, more standard 
iteration schemes (contractions, implicit function theorem, 
schemes based on Duhamel principle, etc.). 
Typical situations where such a loss is unavoidable are related, for example, 
to the presence of the so-called ``small divisors''.
In addition to that, sometimes it could be convenient to use a Nash-Moser iteration 
even if other techniques are also available.
In general, the advantages of the Nash-Moser method for nonlinear PDEs 
(especially quasi-linear ones)
with respect to other approaches are essentially these:
the required estimates on the solution of the linearized problem
allow some loss of regularity, also with respect to the coefficients;
the continuity of the solution of the linearized problem 
with respect to the linearization point 
is not required for the existence proof;
linearizing does not introduce nonlocal terms
(whereas, for example, in some other schemes paralinearizing does);
the nonlinear scheme is ``packaged'' in the theorem and ready-to-use, 
and its application to a PDE problem reduces to verify its assumptions,
which mainly consists of a careful analysis of the linearized operator. 
 
Without claiming to be complete, 
Nash-Moser schemes in Cauchy problems for nonlinear PDEs 
(especially with derivatives in the nonlinearity) 
have been used, for example, 
by Klainerman \cite{Klainerman-1980, Klainerman-1982}
and, more recently, Lindblad \cite{Lindblad}, 
Alvarez-Samaniego and Lannes \cite{Alvarez-Samaniego-Lannes, Lannes}, 
Alexandre, Wang, Xu and Yang \cite{Alexandre}
(see also Mouhot-Villani \cite{Mouhot-Villani}) 
and, in control problems, 
by Beauchard, Coron, Alabau-Boussouira, Olive \cite{Beau1, BC-JFA-2006, Beau2, ACO}
(a discussion about Nash-Moser method in the context of controllability of PDEs 
can be found in \cite{Coron}, section 4.2.2).

The Nash-Moser theorem was first introduced by Nash \cite{Nash}, 
then many refinements, improvements and new versions were developed afterwards:
without demanding completeness we mention, for example, 
the results by Moser \cite{Moser}, Zehnder \cite{Zehnder},
Hamilton \cite{Hamilton}, Gromov \cite{Gromov}, H\"ormander \cite{Geodesy, Olli, H-90}, 
Alinhac and G\'erard \cite{AG}, 
and, more recently, 
Berti, Bolle, Corsi and Procesi \cite{BBPro, BCPro}, 
Texier and Zumbrun \cite{Texier-Zumbrun}, 
Ekeland and S\'er\'e \cite{Ekeland, ES}.

The iteration scheme by H\"ormander \cite{Geodesy} 
(based on telescoping series, and very close to the original scheme by Nash)
is the one used for Cauchy problems by Klainerman \cite{Klainerman-1980, Klainerman-1982}
and by Lindblad \cite{Lindblad}.
\Hor's theorem in \cite{Geodesy} is formulated in the setting of H\oe lder spaces, 
and it also holds for other families of Banach spaces satisfying the same set of basic properties.
Instead, Sobolev spaces do not satisfy that set of properties 
(see Remark \ref{rem:2707.1}).
The same point is expressed, in other words, in \cite{Olli, H-90}.
The theorems in \cite{Olli} and \cite{H-90} are formulated as abstract results, 
with sharp loss of regularity, in the class of \emph{weak} Banach spaces $E_a'$, 
which \Hor{} defines, using smoothing operators, starting from some given scale of Banach spaces $E_a$, $a \geq 0$. 
A key point is that if $E_a$ is a H\oe{}lder space (with exponent $a \notin \N$), 
then it coincides with its weak counterpart $E_a'$, with equivalent norms 
(this is stated explicitly in \cite{Olli}, and proved implicitly in \cite{Geodesy}).
On the contrary, if $E_a$ is a Sobolev space, then $E_a'$ is a strictly larger set, 
with a strictly weaker norm (in fact it is a Besov space, see Remark \ref{rem:2707.2}). 
What is true in Sobolev class is that $E_a \subset E_a' \subset E_b$ for all $b < a$, 
with continuous inclusions.
This is the reason why the application of \Hor's theorems in Sobolev class 
produces a further, unavoidable, arbitrarily small loss. 
This further loss is not present if the theorems of \cite{Olli, H-90}
are applied in the weak spaces $E_a'$,
but these $E_a'$ are not the usual Sobolev spaces 
(see also Remark 1.2 in \cite{BFH}).

In Theorem \ref{thm:1} we overcome this issue by modifying the iteration scheme of \cite{Geodesy},  inspired by a trick commonly used in paradifferential calculus 
(see Remark \ref{rem:11.09.2017}).

\medskip

Theorem \ref{thm:1} is stated in Section \ref{sec:NM}, 
and it is followed by several comments and technical remarks. 
Its proof is contained in Section \ref{sec:proof}. 
An application of the theorem 
is given in Section \ref{sec:appl},
where we remove the loss of regularity from the results in \cite{BFH} 
about control and Cauchy problems for quasi-linear perturbations 
of the Korteweg-de Vries equation in Sobolev class 
(Theorems \ref{thm:control} and \ref{thm:byproduct}).
Possible applications to other PDEs 
are also mentioned (Remark \ref{rem:3108}).

\medskip

\begin{small}
\noindent
\textbf{Acknowledgements}.
This research was supported by the European Research Council under FP7 (ERC Project 306414),
by PRIN 2012 ``Variational and perturbative aspects of nonlinear differential problems'',
and partially by Programme STAR
(UniNA and Compagnia di San Paolo).
\end{small}

\section{A Nash-Moser-\Hor{} theorem} 
\label{sec:NM}

Let $(E_a)_{a \geq 0}$ be a decreasing family of Banach spaces with continuous injections  
$E_b \hookrightarrow E_a$, 
\begin{equation} \label{S0}
\| u \|_a \leq \| u \|_b \quad \text{for} \  a \leq b.	
\end{equation}
Set $E_\infty = \cap_{a\geq 0} E_a$ with the weakest topology making the 
injections $E_\infty \hookrightarrow E_a$ continuous. 
Assume that $S_j : E_0 \to E_\infty$ for $j = 0,1,\ldots$ are linear operators 
such that, with constants $C$ bounded when $a$ and $b$ are bounded, 
and independent of $j$,
\begin{alignat}{2}
\label{S1} 
\| S_j u \|_a 
& \leq C \| u \|_a 
&& \text{for all} \ a;
\\
\label{S2} 
\| S_j u \|_b 
& \leq C 2^{j(b-a)} \| S_j u \|_a 
&& \text{if} \ a<b; 
\\
\label{S3} 
\| u - S_j u \|_b 
& \leq C 2^{-j(a-b)} \| u - S_j u \|_a 
&& \text{if} \ a>b; 
\\ 
\label{S4} 
\| (S_{j+1} - S_j) u \|_b 
& \leq C 2^{j(b-a)} \| (S_{j+1} - S_j) u \|_a 
\quad && \text{for all $a,b$.}
\end{alignat}
From \eqref{S2}-\eqref{S3} one can obtain the logarithmic convexity of the norms
\begin{equation} \label{S5} 
\| u \|_{\lambda a +(1-\lambda) b } 
\leq C\|u\|_a^\lambda\|u\|_b^{1-\lambda} 
\quad \text{if} \ 0 < \lambda < 1.
\end{equation}
Set 
\begin{equation}  \label{new.24}
R_0 u := S_1 u, \qquad 
R_j u := (S_{j+1} - S_j) u, \quad j \geq 1.
\end{equation}
Thus 
\begin{equation} \label{2705.3}
\| R_j u \|_b \leq C 2^{j(b-a)} \| R_j u \|_a \quad \text{for all} \ a,b.
\end{equation}
Bound \eqref{2705.3} for $j \geq 1$ is \eqref{S4}, 
while, for $j=0$, it follows from \eqref{S0} and \eqref{S2}.

We also assume that 
\begin{equation} \label{2705.4}
\| u \|_a^2 \leq C \sum_{j=0}^\infty \| R_j u \|_a^2	\quad \forall a \geq 0,
\end{equation}
with $C$ bounded for $a$ bounded.
This is a sort of ``orthogonality property'' of the smoothing operators.

Now let us suppose that we have another family $F_a$ of decreasing Banach spaces with smoothing operators having the same properties as above. We use the same notation also for the smoothing operators. 

\begin{theorem} \label{thm:1}
Let $a_1, a_2, \a, \b, a_0, \mu$ be real numbers with 
\begin{equation} \label{ineq 2016}
0 \leq a_0 \leq \mu \leq a_1, \qquad 
a_1 + \frac{\b}{2} \, < \a < a_1 + \b , \qquad 
2\a < a_1 + a_2. 
\end{equation}
Let $V$ be a convex neighborhood of $0$ in $E_\mu$. 
Let $\Phi$ be a map from $V$ to $F_0$ such that $\Phi : V \cap E_{a+\mu} \to F_a$ 
is of class $C^2$ for all $a \in [0, a_2 - \mu]$, with 
\begin{align} 
\|\Phi''(u)[v,w] \|_a 
& \leq M_1(a) \big( \| v \|_{a+\mu} \| w \|_{a_0} + \| v \|_{a_0} \| w \|_{a+\mu} \big) 
\notag \\ & \quad 
+ \{ M_2(a) \| u \|_{a+\mu} + M_3(a) \} \| v \|_{a_0} \| w \|_{a_0}
\label{Phi sec}
\end{align}
for all $u \in V \cap E_{a+\mu}$, $v,w \in E_{a+\mu}$,
where $M_i : [0, a_2 - \mu] \to \R$, $i = 1,2,3$, are positive, increasing functions. 
Assume that $\Phi'(v)$, for $v \in E_\infty \cap V$ 
belonging to some ball $\| v \|_{a_1} \leq \d_1$,
has a right inverse $\Psi(v)$ mapping $F_\infty$ to $E_{a_2}$, and that
\begin{equation}  \label{tame in NM}
\| \Psi(v)g \|_a \leq 
L_4(a) \|g\|_{a + \b - \a} + 
\{ L_5(a) \| v \|_{a + \b} + L_6(a) \} \| g \|_0
\quad \forall a \in [a_1, a_2],
\end{equation}
where $L_i : [a_1, a_2] \to \R$, $i = 4,5,6$, 
are positive, increasing functions.

Then for all $A > 0$ there exists $\d > 0$ such that, 
for every $g \in F_\b$ satisfying
\begin{equation} \label{2705.1}
\sum_{j=0}^\infty \| R_j g \|_\b^2 \leq A^2 \| g \|_\b^2, \quad
\| g \|_\b \leq \d,
\end{equation}
there exists $u \in E_\a$ solving $\Phi(u) = \Phi(0) + g$.
The solution $u$ satisfies 
\[
\| u \|_\a \leq C L_{456}(a_2) (1 + A) \| g \|_\b, 
\]
where $L_{456} = L_4 + L_5 + L_6$ 
and $C$ is a constant depending on $a_1, a_2, \a, \b$. 
The constant $\d$ is 
\[
\d = 1/B, \quad 
B = C' L_{456}(a_2) (1+A) \max \big\{ 1, 1/\d_1, L_{456}(a_2) M_{123}(a_2-\mu) \big\}
\]
where $M_{123} = M_1 + M_2 + M_3$ 
and $C'$ is a constant depending on $a_1, a_2, \a, \b$.  

Moreover, let $c > 0$
and assume that \eqref{Phi sec} holds for all $a \in [0, a_2 + c - \mu]$,
$\Psi(v)$ maps $F_\infty$ to $E_{a_2 + c}$, 
and \eqref{tame in NM} holds for all $a \in [a_1, a_2 + c]$. 
If $g$ satisfies \eqref{2705.1} and, in addition, $g \in F_{\b+c}$ with
\begin{equation} \label{0406.1}
\sum_{j=0}^\infty \| R_j g \|_{\b+c}^2 \leq A_c^2 \| g \|_{\b+c}^2 
\end{equation}
for some $A_c$, then the solution $u$ belongs to $E_{\a + c}$, 
with 
\begin{equation} \label{0211.10}	
\| u \|_{\a+c} \leq C_c \big\{ \mG_1 (1+A) \| g \|_\b + \mG_2 (1+A_c) \| g \|_{\b+c} \big\} 
\end{equation}
where 
\begin{align} 
\mG_1 & := \tilde L_6 + \tilde L_{45} (\tilde L_6 \tilde M_{12} + L_{456}(a_2) \tilde M_3) 
\sum_{j=0}^{N-2} z^j, 
\quad \mG_2 := \tilde L_{45} \sum_{j=0}^{N-1} z^j, 
\\
z & := L_{456}(a_1) M_{123}(0) + \tilde L_{45} \tilde M_{12},
\end{align}
$\tilde L_{45} := \tilde L_4 + \tilde L_5$, 
$\tilde L_i := L_i(a_2+c)$, $i = 4,5,6$; 
$\tilde M_{12} := \tilde M_1 + \tilde M_2$, 
$\tilde M_i := M_i(a_2 + c - \mu)$, $i = 1,2,3$;
$N$ is a positive integer depending on $c,a_1,\a,\b$; 
and $C_c$ depends on $a_1, a_2, \a, \b, c$.
\end{theorem}

\subsection{Comments}
\label{subsec:comments}

\begin{remark}
\label{rem:2707.3}
We underline that, in the higher regularity case $g \in F_{\b+c}$, 
the smallness assumption $\| g \|_\b \leq \d$ is only required in ``low'' norm
in Theorem \ref{thm:1} (and $\d$ is independent of $c$).
\end{remark}

\begin{remark}
\label{rem:2707.4}
If the first inequality in \eqref{2705.1} does not hold, 
then one can apply Theorem 2.2.2 in \cite{Geodesy} 
or Theorem 7.1 in \cite{BFH}, obtaining the same type of result with a small additional loss of regularity.
The same if \eqref{2705.4} does not hold.
\end{remark}

\begin{remark}
\label{rem:0911}
With respect to the implicit function theorems in \cite{Geodesy,Olli,BFH}, in Theorem \ref{thm:1} we slightly modify the form of the tame estimates concerning $\Phi''$ and $\Psi$, allowing the presence of extra terms, corresponding to $M_3(a)$ in \eqref{Phi sec} and $L_6(a)$ in \eqref{tame in NM}. The introduction of these terms is natural when one is interested in keeping explicitly track of the high operator norms of $\Phi$.
\end{remark}

\begin{remark}
\label{rem:param}
Theorem \ref{thm:1} could also be stated with $a_0 = \mu = a_1$, 
since in the proof $a_0, \mu$ are often deteriorated to $a_1$. 
However, in the applications to PDEs, $\Phi$ is usually a differential operator, 
and in principle it is somewhat natural to distinguish its loss of regularity $\mu$ 
(the order of $\Phi$), 
the low norm threshold $a_0$ appearing in the tame estimates \eqref{Phi sec}
(usually given by the $L^\infty$ embedding),
and the minimal regularity $a_1$ at which the linearized operator $\Phi'(v)$ 
admits a right inverse $\Psi(v)$. 
Regarding the other parameters of the theorem, $a_2$ is the ``high'' norm 
required by the proof of the first part of the theorem, giving the solution $u \in E_\a$; 
$\b - \a$ is the loss of regularity of $\Psi(v)h$ in terms of its argument $h$
(namely the order of the operator $\Psi(v)$),
and $\b$ is the loss of regularity of $\Psi(v)$ in terms of its coefficient $v$
($v$ is the point where $\Phi$ has been linearized), see \eqref{tame in NM}. 
Thus in the thesis of Theorem \ref{thm:1} 
$\b$ is the regularity of the datum $g$, 
and $\a$ is the one of the solution $u$ of the equation $\Phi(u) = \Phi(0) + g$.

Note that, given $g \in F_\b$, and given $v \in E_\infty$ with $\| v \|_{a_1} \leq \d_1$,
the linearized equation $\Phi'(v) h = g$ has a solution $h = \Psi(v) g \in E_\a$ 
(see \eqref{tame in NM}); 
hence the solution $u \in E_\a$ of the \emph{nonlinear} equation $\Phi(u) = \Phi(0) + g$ 
given by Theorem \ref{thm:1} has the same regularity 
as the solution of the linearized problem with the same datum.
In this sense our theorem is sharp: the nonlinear problem reaches exactly 
the same regularity given by the linearized one. 
\end{remark}

\begin{remark} 
\label{rem:2707.2}
As already said in the Introduction, 
if $E_a$ is a Sobolev space $H^a$, 
then the weak space $E_a'$ defined in \cite{Olli}
is a strictly larger set, with a strictly weaker norm, and it is in fact the Besov space
$B^a_{2,\infty}$.
To show it, we start by recalling the general definition of $E_a'$ in \cite{Olli}. 

\medskip

\noindent
\emph{Definition of $E_a'$ in \cite{Olli}}.
Assume that $(E_a)_{a \geq 0}$ is a family of Banach spaces, 
with $E_b \subset E_a$, $\| u \|_a \leq \| u \|_b$ for $a < b$. 
Let $E_\infty = \cap_{a \geq 0} E_a$. 
Assume that $S_\theta : E_0 \to E_\infty$, with real parameter $\theta \geq 1$, 
is a family of linear operators such that, with constants $C$ bounded for $a,b$ bounded, 

$(i)$ $\| S_\theta u \|_b \leq C \| u \|_a$ for $b \leq a$;

$(ii)$ $\| S_\theta u \|_b \leq C \theta^{b-a} \| u \|_a$ for $a < b$;

$(iii)$ $\| u - S_\theta u \|_b \leq C \theta^{b-a} \| u \|_a$ for $a > b$;

$(iv)$ $\| \frac{d}{d\theta} S_\theta u \|_b \leq C \theta^{b-a-1} \| u \|_a$ for all $a,b$.

\noindent
Consider an increasing sequence $1 = \theta_0 < \theta_1 < \ldots \to \infty$ 
with $\theta_{j+1} / \theta_j$ bounded, 
and let $\Delta_j = \theta_{j+1} - \theta_j$.  
Let $a_1 < a < a_2$. 
Then $E_a'$ is defined in \cite{Olli} as 
the set of all sums $u = \sum_{j=0}^\infty \Delta_j u_j$, 
with $u_j \in E_{a_2}$, for which there exists $M > 0$ 
such that, for all $j \in \N$, 
\[
\| u_j \|_{a_1} \leq M \theta_j^{a_1 - a - 1}, 
\quad
\| u_j \|_{a_2} \leq M \theta_j^{a_2 - a - 1}.
\]
The norm $\| u \|_{E_a'}$ is defined in \cite{Olli} as 
the infimum of $M$ over all such decompositions.

\medskip

In \cite{Olli} it is also observed 
that, up to equivalent norms, it is sufficient to calculate $M$ 
for the decomposition defined by $u_j = R_j u$, where $R_0 u = S_{\theta_1} u / \Delta_0$ 
and $R_j u = (S_{\theta_{j+1}} u - S_{\theta_j} u) / \Delta_j$ for $j \geq 1$;
that $E_a \subset E_a' \subset E_b$ for all $b < a$,  
with continuous inclusions;
that different choices of the family $S_{\theta_j}$ 
lead to the same set $E_a'$ with equivalent norms;
that different choices of $a_1, a_2$ with $a_1 < a < a_2$ 
also lead to the same set $E_a'$ with equivalent norms.

When $(E_a)$ is the family of Sobolev spaces on $\R^d$
\[
E_a = H^a(\R^d,\C) := \Big\{ 
u(x) = \int_{\R^d} \hat u(\xi) \, e^{i \xi \cdot x} \, d\xi : 
\| u \|_a^2 := \int_{\R^d} |\hat u(\xi)|^2 \, \langle \xi \rangle^{2a} \, d\xi 
< \infty \Big\},
\]
where $\langle \xi \rangle := (1 + |\xi|^2)^{\frac12}$, 
or on $\T^d$
\[
E_a = H^a(\T^d,\C) := \Big\{ 
u(x) = \sum_{k \in \Z^d} \hat u_k \, e^{ik \cdot x} : 
\| u \|_a^2 := \sum_{k \in \Z^d} |\hat u_k|^2 \, \langle k \rangle^{2a} < \infty \Big\},
\]
where $\T := \R / 2 \pi \Z$, 
one can define $S_\theta$ as the smooth Fourier cut-off operator
\[
S_\theta u(x) = \int_{\R^d} \hat u(\xi) \psi\Big( \frac{|\xi|}{\theta} \Big) e^{i \xi \cdot x} \, d\xi
\quad \text{or} \quad 
S_\theta u(x) = \sum_{k \in \Z^d} \hat u_k \psi\Big( \frac{|k|}{\theta} \Big) e^{ik\cdot x},
\]
where $\psi \in C^\infty$, $0 \leq \psi \leq 1$, 
$\psi = 1$ on $[0,1]$ and $\psi = 0$ on $[2,\infty)$.
One can easily check that properties $(i),(ii),(iii),(iv)$ are satisfied. 
Then, taking $\theta_j = 2^j$, the sum $u = \sum_{j = 0}^\infty \Delta_j R_j u$
defined above is a Littlewood-Paley decomposition of $u$. 
It follows that $\| u \|_{E_a'}$ is equivalent to $\sup_{j \geq 0} \| \Delta_j R_j u \|_{a}$,
which is the $\ell^\infty$ norm of the sequence of the Sobolev norms 
of the dyadic blocks of $u$,
so that $E_a'$ is the Besov space $B^a_{p,r}$ with $p = 2$ and $r = \infty$.
Since $\| u \|_a$ is equivalent to the $\ell^2$ norm of the same sequence, 
and $\| \ \|_{\ell^\infty} \leq \| \ \|_{\ell^2}$, 
it follows that the norm of $E_a'$ is weaker than the one of $E_a$.
Moreover $E_a$ is strictly contained in $E_a'$ because $\ell^2$ is strictly contained in $\ell^\infty$:
for example, the function 
\begin{equation} \label{2807.1}
u(x) := \int_{\R^d} \langle \xi \rangle^{-a - \frac{d}{2}} \, e^{i \xi \cdot x} \, d\xi
\quad \text{or} \quad 
u(x) := \sum_{k \in \Z^d} \langle k \rangle^{-a - \frac{d}{2}} \, e^{ik \cdot x}
\end{equation}
belongs to $E_a' \setminus E_a$, because the sequence of the $H^a$ norm of its dyadic blocks
is in $\ell^\infty \setminus \ell^2$, as one can check directly.
\end{remark}

\begin{remark} \label{rem:2707.1}
In Appendix A of \cite{Geodesy}, H\oe{}rmander discusses various properties of families 
of H\oe{}lder spaces $C^{k,\a}(B)$ where $B$ is a compact convex subset of $\R^n$
with nonempty interior. 
Among other results, it is shown in \cite{Geodesy} that the spaces $\mH^a$ 
with real parameter $a \geq 0$, defined by $\mH^0 := C(B)$ and 
$\mH^a := C^{k,\a}(B)$ with $k + \a = a$, $0 < \a \leq 1$, and $k \geq 0$ integer,
form a family of Banach spaces to which H\oe{}rmander's implicit function Theorem 2.2.2 of \cite{Geodesy} applies. 
On the contrary, some of the key results of Appendix A of \cite{Geodesy}
do not hold for families of Sobolev spaces.
In particular, this is the case for Theorem A.11 in \cite{Geodesy}, 
which is stated for $\mH^a = C^{k,\a}(B)$ in the case $0 < \a < 1$:

\medskip

\noindent
\textbf{Theorem A.11 of \cite{Geodesy}}.
\emph{Let $u_\theta$ for $\theta > \theta_0$ be a $C^\infty$ function in $B$
and assume that $\| u_\theta \|_{a_i} \leq M \theta^{b_i - 1}$, $i=0, 1$,
where  $b_0 < 0 < b_1$ and $a_0 < a_1$. 
Define $\lm$ by $\lm b_0 + (1-\lm) b_1 = 0$ 
and set $a = \lm a_0 + (1-\lm) a_1$,
that is,
$a = (a_0 b_1 - a_1 b_0)/(b_1 - b_0)$.
If $a = k + \a$ with $k$ integer and $0 < \a < 1$ 
(so that $a$ is not an integer), 
it follows then that $u = \int_{\theta_0}^\infty u_\theta \, d\theta$ 
is in $\mH^a = C^{k,\a}(B)$ and $\| u \|_a = \| u \|_{C^{k,\a}(B)} \leq C_a M$.}

\medskip

It is not difficult to see that a corresponding result for Sobolev spaces does not hold.
For example, in the Sobolev space $H^s(\R^d,\C)$ take
$u_\theta(x) = \int_{\R^d} \ph(|\xi| / \theta) \, e^{i \xi \cdot x} \, d\xi \, \theta^{-\beta}$ 
where 
$\ph \in C^\infty(\R)$, supp$(\ph) \subseteq [\frac12, \frac32]$, 
with $0 \leq \ph \leq 1$, and $\ph(1) = 1$. 
Let $\b > \frac{d}{2} + 1$, $\theta_0 = 1$, and fix $a_0, a_1$ such that
$0 \leq a_0 < \b - \frac{d}{2} - 1 < a_1$.
Let $b_i := a_i - \b + \frac{d}{2} + 1$, $i = 0,1$,
so that $b_0 < 0 < b_1$. 
Hence $u_\theta$ satisfies the estimates $\| u_\theta \|_{a_i} \leq M \theta^{b_i - 1}$,
and $a := (a_0 b_1 - a_1 b_0) / (b_1 - b_0)$ is given by $a = \b - \frac{d}{2} - 1$. 
However, the function $u = \int_1^\infty u_\theta \, d\theta$
has Fourier transform $\hat u(\xi) = \int_1^\infty \ph(|\xi| / \theta) \theta^{-\b} \, d\theta$.
Now $| \hat u(\xi)| \geq C |\xi|^{1-\b}$ for all $|\xi| \geq 1$, 
and therefore $| \hat u(\xi)| |\xi|^a \geq C |\xi|^{-\frac{d}{2}}$, 
whence $u \notin H^a(\R^d,\C)$. 

Similarly, on the Sobolev space $H^s(\T^d, \C)$ of periodic functions, 
we take 
\[
u_\theta(x) = \sum_{k \in \Z^d, \frac12 \theta \leq |k| \leq \frac32 \theta}
e^{i k \cdot x} \, \theta^{-\beta}.
\]
Let $\b, a_0, a_1, b_0, b_1, a$ as above. 
Hence $u_\theta$ satisfies the estimates $\| u_\theta \|_{a_i} \leq M \theta^{b_i - 1}$.
The function $u = \int_1^\infty u_\theta \, d\theta$
has Fourier coefficients $\hat u_k = \int_{\frac23 |k|}^{2|k|} \theta^{-\b} \, d\theta
\geq C |k|^{1-\b}$. 
Therefore $|\hat u_k| |k|^a \geq C |k|^{-\frac{d}{2}}$, 
whence $u \notin H^a(\T^d, \C)$. 

A consequence, Theorem A.11 of \cite{Geodesy} does not hold for Sobolev spaces
and hence Theorem 2.2.2 of \cite{Geodesy} does not apply.
\end{remark}

\begin{remark} \label{rem:2107.7} 
We make an attempt to discuss the consequences of the ``velocity'' of the sequence 
$(\theta_j)$ of smoothing operators in different Nash-Moser theorems.

In Moser \cite{Moser-Pisa-66-parte-1}, Zehnder \cite{Zehnder},
and recent improvements like \cite{BBPro, ES}, 
the sequence $S_{\theta_j}$ of smoothing operators along the iteration scheme 
is defined as $\theta_{j+1} = \theta_j^\chi$, with $1 < \chi < 2$
($\chi = \frac32$ in \cite{Moser-Pisa-66-parte-1}), 
namely 
\[
\theta_j = \theta_0^{\chi^j}
\]
with $\theta_0 > 1$.
Thus $\theta_j$, the ratio $\theta_{j+1} / \theta_j$ and the difference $\theta_{j+1} - \theta_j$
all diverge to $\infty$ as $j \to \infty$. 

On the opposite side, in \Hor{} \cite{Geodesy, Olli, H-90} the ``velocity'' 
of the smoothings is 
\[
\theta_j = (a + j)^\e
\] 
with $a > 0$ large and $\e \in (0,1)$ small, 
so that $\theta_j$ diverges, 
the ratio $\theta_{j+1} / \theta_j$ tends to 1 
and the difference $\theta_{j+1} - \theta_j$ goes to zero.
This choice corresponds to a very fine discretization 
of the continuous real parameter $\theta \in [1, \infty)$ of Nash \cite{Nash}.

An intermediate choice is 
\[
\theta_j = c^j
\]
for some $c > 1$.
In this case $\theta_j \to \infty$, 
the ratio $\theta_{j+1} / \theta_j$ is constant and equal to $c$, 
and the difference $\theta_{j+1} - \theta_j \to \infty$.
This is the choice in \cite{Klainerman-1982} with $c = 2^\e$ 
(equations (4.4), (S1), (S2) in \cite{Klainerman-1982}).
For $c = 2$, it corresponds to the dyadic Littlewood-Paley decomposition, 
and it is our choice in Theorem \ref{thm:1}. 

\medskip

The velocity of the sequence $\theta_j$ 
has the following consequences.

\medskip

If the ratio $\theta_{j+1} / \theta_j$ diverges to $\infty$, 
then a further loss of regularity is \emph{introduced} in the process of constructing the solution.
The main reason of this artificial loss is that the high and low norms of the difference 
$(S_{\theta_{j+1}} - S_{\theta_j}) u$ cannot be sharply estimated 
in terms of the corresponding powers of $\theta_j$ only, 
but, instead, one has 
\begin{equation} \label{2807.2}
\frac{1}{\theta_{j+1} - \theta_j}\, \| (S_{\theta_{j+1}} - S_{\theta_j}) u \|_b 
\leq C_{a,b} \max \{ \theta_j^{b-a-1} , \theta_{j+1}^{b-a-1} \} \| u \|_a ,
\end{equation}
and the maximum is $\theta_j^{b-a-1}$ or $\theta_{j+1}^{b-a-1}$ 
according to the (high or low) norm one is estimating.
Along the iteration scheme one has to estimate both high and low norms, 
and the discrepancy between $\theta_{j+1}^{b-a-1}$ and $\theta_j^{b-a-1}$
generates a loss of regularity. 
In the particular case $\theta_{j+1} = \theta_j^\chi$, for $b > a +1$ 
one can write \eqref{2807.2} in terms of an explicit loss $\s$ of regularity, namely
\begin{equation} \label{2807.3}
\frac{1}{\theta_{j+1} - \theta_j}\, \| (S_{\theta_{j+1}} - S_{\theta_j}) u \|_b 
\leq C_{a,b} \theta_j^{b-a-1+\s} \| u \|_a
\end{equation}
where $(\chi-1)(b-a-1) \leq \s$.

Instead, when the ratio $\theta_{j+1} / \theta_j$ is bounded, 
\eqref{2807.2} reduces to 
\begin{equation} \label{2807.4}
\frac{1}{\theta_{j+1} - \theta_j}\, \| (S_{\theta_{j+1}} - S_{\theta_j}) u \|_b 
\leq C_{a,b} \theta_j^{b-a-1} \| u \|_a.
\end{equation}

\medskip

If the difference $\theta_{j+1} - \theta_j$ tends to zero, 
then this can be used to simplify the proof of the convergence 
of the quadratic error in the telescoping \Hor{} scheme.
This is what is done in \cite{Olli} to obtain bound (15), 
where $\theta_{j+1} - \theta_j = O(\theta_j^{-N})$
and $N$ is chosen large enough.
See also the estimate of the term $e_k''$ on page 150 in Alinhac-G\'erard \cite{AG}.

\medskip

In Sobolev class, 
the orthogonality property \eqref{2705.4} is somehow related to the velocity of $\theta_j$ 
in the following sense. 
Consider $E_a = H^a(\T^d)$ or $H^a(\R^d)$. 
If $S_\theta$ is the ``crude'' Fourier truncation operator
\[
S_\theta u(x) = \sum_{k \in \Z^d, |k| \leq \theta} \hat u_k e^{ik\cdot x}
\quad \text{or} \quad 
S_\theta u(x) = \int_{|\xi| \leq \theta} \hat u(\xi) e^{i \xi \cdot x} \, d\xi,
\]
and $R_0 := S_{\theta_1}$, $R_j := (S_{\theta_{j+1}} - S_{\theta_j})$, 
then \eqref{2705.4} holds no matter what the choice of the sequence $\theta_j$ is 
(with $\theta_0 < \theta_1 < \theta_2 < \ldots \to \infty$).

If, instead, $S_\theta$ is a smooth Fourier cut-off operator
\[
S_\theta u(x) = \sum_{k \in \Z^d} \hat u_k \psi\Big( \frac{|k|}{\theta} \Big) e^{ik\cdot x}
\quad \text{or} \quad 
S_\theta u(x) = \int_{\R^d} \hat u(\xi) \psi\Big( \frac{|\xi|}{\theta} \Big) e^{i \xi \cdot x} \, d\xi,
\]
where $\psi \in C^\infty$, $0 \leq \psi \leq 1$, 
$\psi = 1$ on $[0,1]$ and $\psi = 0$ on $[2,\infty)$,
then the orthogonality condition \eqref{2705.4} holds if 
$\theta_{j+1} / \theta_j \geq c > 1$, 
and it does not hold if $\theta_{j+1} / \theta_j \to 1$.
These smooth Fourier cut-offs, commonly used in Fourier analysis,
are a natural choice when property $(iv)$ of \cite{Olli} has to be satisfied
(properties $(i)$-$(iv)$ of \cite{Olli} are recalled in Remark \ref{rem:2707.2};
in Theorem \ref{thm:1}, property $(iv)$ of \cite{Olli} has been replaced 
by the less demanding inequality \eqref{S4}).
\end{remark}

\section{Proof of Theorem \ref{thm:1}}
\label{sec:proof}

Fix $\g > 0$ such that $2 a_1 + \b + \g \leq 2\a$.    
In this proof we denote by $C$ any constant (possibly different from line to line) 
depending only on $a_1, a_2, \a, \b, \mu, a_0, \g$, which are fixed parameters. 
Denote, in short, 
\begin{equation} \label{3005.4}
g_j := R_j g \quad \forall j \geq 0. 
\end{equation}
By \eqref{2705.3}, 
\begin{equation} \label{2705.2} 
\| g_j \|_b \leq C_b \, 2^{j(b-\b)} \| g_j \|_\b \quad \forall b \in [0, +\infty).
\end{equation}

\noindent
\textbf{Recursive scheme.}
We claim that, if $\| g \|_\b$ is small enough, 
then we can define a sequence $u_j \in V \cap E_{a_2 + c}$ with $u_0 := 0$ 
by the recursion formula
\begin{equation} \label{Olli.8}
u_{j+1} := u_j + h_j, \quad
v_j := S_j u_j, \quad
h_j := \Psi(v_j) (g_j + y_j) \quad \forall j \geq 0,
\end{equation}
where $y_0 := 0$, 
\begin{equation} \label{new.1}
y_1 := - S_1 e_0, \qquad
y_j := - S_j e_{j-1} - R_{j-1} \sum_{i=0}^{j-2} e_i 
\quad \ \forall j \geq 2,
\end{equation}
and $e_j := e_j' + e_j''$, 
\begin{equation} \label{new.2}
e_j' := \Phi(u_j + h_j) - \Phi(u_j) - \Phi'(u_j) h_j , 
\qquad
e_j'' := (\Phi'(u_j) - \Phi'(v_j)) h_j.
\end{equation}
The fact that the recursive scheme \eqref{Olli.8}-\eqref{new.2} is well-defined
will be a consequence of the following estimates.

\medskip

\noindent
\textbf{Iterative estimates.} 
We prove that there exist positive constants $K_1, \ldots, K_4$ 
such that, for all $j \geq 0$,
\begin{align} 
\| h_j \|_a 
& \leq K_1 (\| g \|_\b \, 2^{-j \g} + \| g_j \|_\b) \, 2^{j(a-\a)} \quad 
\forall a \in [a_1, a_2], 
\label{Olli.9}
\vspace{2pt} \\
\| v_j \|_a 
& \leq K_2 \| g \|_\b \, 2^{j(a-\a)} \quad 
\forall a \in [a_1 + \b, a_2 + \b], 
\label{Olli.10}
\vspace{2pt} \\
\| u_j - v_j \|_a 
& \leq K_3 \| g \|_\b \, 2^{j(a-\a)} \quad
\forall a \in [0, a_2],
\label{Olli.11}
\vspace{2pt} \\
\| u_j \|_\a 
& \leq K_4 \| g \|_\b.
\label{reve}
\end{align} 
We prove \eqref{Olli.9}-\eqref{reve} by induction. 

\medskip

\noindent
\textsc{Base case.}
For $j=0$, \eqref{Olli.10}, \eqref{Olli.11} and \eqref{reve} are trivially satisfied, 
and \eqref{Olli.9} follows from \eqref{2705.2} because $h_0 = \Psi(0)g_0$, 
provided that $C (L_4(a_2) + L_6(a_2)) \leq K_1$.

\medskip

\noindent
\textsc{Inductive step.}
Let $k \geq 0$ and assume that, for all $j=0,\ldots,k$, 
\eqref{Olli.9}, \eqref{Olli.10}, \eqref{Olli.11}, \eqref{reve} hold.  

\noindent
$\bullet$ \emph{Proof of \eqref{reve} at $j=k+1$.} 
By \eqref{2705.3} and \eqref{Olli.9} one has for all $n \leq k$, all $j \geq 0$,
\begin{equation} \label{3005.1}
\| R_j h_n \|_\a 
\leq C \, 2^{j(\a - a)} \| h_n \|_a 
\leq C K_1 \xi_n \, 2^{(j-n)(\a - a)} 
\quad \forall a \in [a_1, a_2],
\end{equation}
where $\xi_n := \| g \|_\b \, 2^{-n \g} + \| g_n \|_\b$.
Since $u_{k+1} = \sum_{n=0}^k h_n$, 
using \eqref{3005.1} with $a = a_1$ if $n > j$ 
and $a = a_2$ if $n \leq j$, we get
\begin{equation} \label{3005.2}
\| R_j u_{k+1} \|_\a 
\leq \sum_{n=0}^k \| R_j h_n \|_\a 
\leq C K_1 (\e_j' + \e_j'')
\end{equation}
where
\begin{equation} \label{3005.3}
\e_j' := \sum_{n = j+1}^k \xi_n \, 2^{-(n-j)(\a - a_1)}, 
\quad  
\e_j'' := \sum_{n = 0}^{ \min\{ k,j \} } \xi_n \, 2^{-(j-n)(a_2 - \a)}
\end{equation}
and $\e_j' = 0$ for $j+1 > k$ (empty sum). 
By H\"older inequality, 
\begin{align}
\sum_{j=0}^\infty \e_j'^2 
& \leq \sum_{j=0}^\infty \Big( \sum_{n=j+1}^k \xi_n^2 \, 2^{-(n-j)(\a - a_1)} \Big)
\Big( \sum_{n=j+1}^k 2^{-(n-j)(\a - a_1)} \Big)
\notag 
\\
& \leq C \sum_{j=0}^\infty \sum_{n=j+1}^k \xi_n^2 \, 2^{-(n-j)(\a - a_1)}
\, = C \sum_{n=1}^k \xi_n^2 \, \sum_{j=0}^{n-1} 2^{-(n-j)(\a - a_1)}
\notag 
\\ 
& \leq C \sum_{n=1}^k \xi_n^2 
\leq C (1 + A)^2 \| g \|_\b^2,
\label{0106.1}
\end{align}
where the last inequality follows from \eqref{2705.1}, \eqref{0406.1} and \eqref{3005.4}.
Similarly, one proves that 
$\sum_{j=0}^\infty \e_j''^2 \leq C (1 + A)^2 \| g \|_\b^2$.
Thus by \eqref{2705.4} and \eqref{3005.2} we deduce that 
\begin{equation} \label{3005.5}
\| u_{k+1} \|_\a \leq C K_1 (1 + A) \| g \|_\b,
\end{equation}
which gives \eqref{reve} if $C K_1 (1 + A) \leq K_4$.

\noindent
$\bullet$ \emph{Proof of \eqref{Olli.11} at $j=k+1$.} 
By \eqref{S3}, \eqref{S1} and \eqref{3005.5} one has
\begin{equation} \label{3005.6}
\| u_{k+1} - v_{k+1} \|_0 
\leq C \, 2^{-(k+1)\a} \| u_{k+1} \|_\a 
\leq C K_1 (1 + A) \| g \|_\b \, 2^{-(k+1)\a}.
\end{equation}
By triangular inequality, \eqref{S1} and \eqref{Olli.9} we get
\begin{align} 
\| u_{k+1} - v_{k+1} \|_{a_2}
& \leq C \| u_{k+1} \|_{a_2} 
\leq C \sum_{n=0}^k \| h_n \|_{a_2} 
\leq C K_1 \| g \|_\b \, 2^{(k+1)(a_2 - \a)}.
\label{3005.7}
\end{align}
Interpolating between $0$ and $a_2$ by \eqref{S5} gives 
$\| u_{k+1} - v_{k+1} \|_a \leq C K_1 (1 + A) \| g \|_\b \, 2^{(k+1)(a - \a)}$ 
for all $a \in [0,a_2]$. 
This gives \eqref{Olli.11} if $C K_1 (1 + A) \leq K_3$.

\noindent
$\bullet$ \emph{Proof of \eqref{Olli.10} at $j=k+1$.}  
We use the assumption $a_1 + \b > \a$, \eqref{S2} and \eqref{3005.5} and we get
\[
\| v_{k+1} \|_a 
\leq C \, 2^{(k+1)(a - \a)} \| u_{k+1} \|_\a
\leq C K_1 (1 + A) \| g \|_\b \, 2^{(k+1)(a -\a)}
\]
for all $a \in [a_1 + \b, a_2 + \b]$. 
This gives \eqref{Olli.10} if $C K_1 (1 + A) \leq K_2$.

\noindent
$\bullet$ \emph{Proof of \eqref{Olli.9} at $j=k+1$.} 
We begin with proving the following estimate of $y_{k+1}$.

\bigskip

\noindent
\textbf{Claim.} \emph{One has}
\begin{equation} \label{new.8}
\| y_{k+1} \|_b 
\leq C K_1 (K_1 + K_3) M_{123}(a_2 - \mu) \| g \|_\b^2 \, 2^{(k+1)(b-\b-\g)} 
\quad \forall b \in [0,a_2 + \b - \a].
\end{equation}
\emph{Proof of Claim \eqref{new.8}.}
Since $u_j, v_j, u_j + h_j$ belong to $V$ for all $j = 0,\ldots, k$, 
we use Taylor formula and \eqref{Phi sec} to deduce that, for $j = 0, \ldots, k$ 
and $a \in [0, a_2 - \mu]$,
\begin{align} \label{new.9}
\| e_j \|_a 
& \leq \| h_j \|_{a+\mu} \| h_j \|_{a_0} \{ M_1(a) + M_2(a) \| h_j \|_{a_0} \}
+ \| h_j \|_{a_0}^2 \{ M_3(a) + M_2(a) \| u_j \|_{a+\mu} \} 
\notag \\ & \quad \ 
+ \| h_j \|_{a_0} \| v_j - u_j \|_{a+\mu} \{ M_1(a) + M_2(a) \| v_j - u_j \|_{a_0} \}
+ \| h_j \|_{a + \mu} \| v_j - u_j \|_{a_0} M_1(a)
\notag \\ & \quad \ 
+ \| h_j \|_{a_0} \| v_j - u_j \|_{a_0} \{ M_3(a) + M_2(a) \| v_j \|_{a+\mu} \}. 
\end{align}
Let $p := \max \{ 0, \b - \a + \mu \}$.  
For future convenience, note that $p \leq a_1 + \b - \a$ because 
$0 < a_1 + \b - \a$ and $\mu + \b - \a \leq a_1 + \b - \a$. 
By assumption, $\g \leq 2\a - \b - 2a_1$ and $2\a - a_1 < a_2$. 
Hence
\begin{equation} \label{0506.2}
\a + p + \g \leq 3\a + p - \b - 2a_1
\leq 3\a + (a_1 + \b - \a) - \b -2a_1 
= 2\a - a_1 < a_2.
\end{equation}
Let $q := a_2 + \b - \a + \mu - p$ 
(so that $q = a_2$ if $\b - \a + \mu \geq 0$, and $q < a_2$ if $\b - \a + \mu < 0$).
For $j = 1, \ldots, k$, by \eqref{Olli.9} we have 
\begin{equation} \label{new.11}
\| u_j \|_q \leq \| u_j \|_{a_2} 
\leq \sum_{i=0}^{j-1} \| h_i \|_{a_2}
\leq K_1 \| g \|_\b \sum_{i=0}^{j-1} 2^{i(a_2-\a)}
\leq C K_1 \| g \|_\b \, 2^{j(a_2-\a)},
\end{equation}
while for $j=0$ we have $u_0 = 0$ by assumption. 
We consider \eqref{new.9} with $a = q - \mu$ 
(note that $q - \mu \in [0, a_2 - \mu]$). 
Since $a_0 \leq a_1$, using \eqref{new.11}, \eqref{Olli.9}, \eqref{Olli.11}
we have 
\begin{align*} 
\| e_j \|_{a_2 + \b - \a - p}
& \leq C K_1 (K_1 + K_3) \| g \|_\b^2 \,
\Big\{ M_1 (a_2 - \mu) 2^{j(a_1 + q - 2\a)} 
\\ & \quad \  
+ M_2 (a_2 - \mu) 2^{j(a_2 + 2a_1 - 3\a)} 
+ M_3 (a_2 - \mu) 2^{j (2 a_1 - 2\a)} \Big\}
\end{align*} 
provided that $K_1 \| g \|_\b \leq 1$. 
We assume that $K_1 \| g \|_\b \leq 1$.
By the definition of $q$, the exponents 
$(a_1 + q - 2\a)$, $(a_2 + 2a_1 - 3\a)$ and $(2 a_1 - 2\a)$
are $\leq (a_2 - \a - p - \g)$ because, by assumption, 
$2a_1 + \b + \g \leq 2\a$.
Thus 
\begin{equation} \label{0506.1}
\| e_j \|_{a_2 + \b - \a - p}
\leq C K_1 (K_1 + K_3) M_{123}(a_2 - \mu) \| g \|_\b^2 \, 
2^{j(a_2 - \a - p - \g)}.
\end{equation}
Now we estimate $\| S_{k+1} e_k \|_0$.
By \eqref{reve}, $\| u_k \|_{\mu} \leq \| u_k \|_\a \leq K_4 \| g \|_\b$, 
and we assume that $K_4 \| g \|_\b \leq 1$.
Since $a_0, \mu \leq a_1$, by \eqref{S1}, \eqref{Olli.9}, \eqref{Olli.11} and \eqref{new.9}, 
using the bound $2a_1 + \b + \g \leq 2\a$, we get
\begin{equation} \label{new.16}
\| S_{k+1} e_k \|_0 
\leq C K_1 (K_1 + K_3) M_{123}(0) \| g \|_\b^2 \, 2^{-(k+1)(\b+\g)}. 
\end{equation}
By \eqref{S2} and \eqref{new.16} we deduce that  
\begin{equation} \label{new.17}
\| S_{k+1} e_k \|_b 
\leq 
C K_1 (K_1 + K_3) M_{123}(0) \| g \|_\b^2 \, 2^{(k+1)(b-\b-\g)} 
\end{equation}
for all $b \in [0, a_2 + \b - \a]$. 
Now we estimate the other terms in $y_{k+1}$ (see \eqref{new.1}). 
For all $b \in [0, a_2 + \b - \a]$, by \eqref{2705.3} and \eqref{0506.1} we have 
\begin{align} 
& \sum_{i=0}^{k-1} \| R_k e_i \|_b
\leq \sum_{i=0}^{k-1} C \, 2^{k(b - a_2 - \b + \a + p)} \| e_i \|_{a_2 + \b - \a - p}
\notag \\ 
& \qquad 
\leq C K_1 (K_1 + K_3) M_{123}(a_2 - \mu) \| g \|_\b^2 \,
2^{k(b - a_2 - \b + \a + p)} 
\sum_{i=0}^{k-1} 2^{i(a_2 - \a - p - \g)}
\notag \\ 
& \qquad 
\leq C K_1 (K_1 + K_3) M_{123}(a_2 - \mu) \| g \|_\b^2 \,
2^{k(b - \b - \g)}
\label{new.18}
\end{align}
because $a_2 - \a - p - \g > 0$ (see \eqref{0506.2}). 
The sum of \eqref{new.17} and \eqref{new.18} completes the proof of Claim \eqref{new.8}.

\bigskip

Now we are ready to prove \eqref{Olli.9} at $j=k+1$. 
By \eqref{S1} and \eqref{3005.5} we have
$\| v_{k+1} \|_{a_1} \leq C \| u_{k+1} \|_{a_1} \leq C K_1 (1 + A) \| g \|_\b$,
and we assume that $C K_1 (1 + A) \| g \|_\b \leq \d_1$, so that $\Psi(v_{k+1})$ is defined.
By \eqref{Olli.8}, \eqref{tame in NM}, \eqref{2705.2}, \eqref{new.8}, \eqref{Olli.10}
one has, for all $a \in [a_1, a_2]$,
\begin{align} 
\| h_{k+1} \|_a 
& \leq C \big\{ K_1 (K_1 + K_3) M_{123}(a_2-\mu) \| g \|_\b^2 \, 2^{-(k+1)\g} + \| g_{k+1} \|_\b \big\} \, 
\notag \\ & \quad 
\cdot \big\{ [L_4(a) + L_5(a)] 2^{(k+1)(a-\a)} + L_6(a) 2^{-(k+1)\b} \big\}
\label{new.21}
\end{align}
if $K_2 \| g \|_\b \leq 1$. 
We assume that $K_2 \| g \|_\b \leq 1$. 
Since $-\b < a_1 - \a$, 
bound \eqref{new.21} implies \eqref{Olli.9} if 
\[
C L_{456}(a_2) \leq K_1, \quad 
C L_{456}(a_2) (K_1 + K_3) M_{123}(a_2 - \mu) \| g \|_\b \leq 1.
\]

\noindent
$\bullet$ \emph{Choice of the constants.}
The induction proof of \eqref{Olli.9}, \eqref{Olli.10}, \eqref{Olli.11}, \eqref{reve} 
is complete if $K_1, K_2, K_3, K_4, \| g \|_\b$ satisfy: 
\begin{align}
& C_* L_{456}(a_2) \leq K_1; \quad  
C_* K_1 (1+A) \leq K_i \ \  \text{for} \ i=2,3,4; \quad 
K_m \| g \|_\b \leq 1 \ \  \text{for} \ m=1,2,4;
\notag \\ 
& C_* K_1 (1 + A) \| g \|_\b \leq \d_1; \quad 
C_* M_{123}(a_2 - \mu) L_{456}(a_2) (K_1 + K_3) \| g \|_\b \leq 1
\label{3110.1}
\end{align}
where $C_*$ is the largest of the constants appearing above. 
First we fix $K_1 = C_* L_{456}(a_2)$. 
Then we fix $K_2 = K_3 = K_4 = C_* K_1 (1+A)$, 
and finally we fix $\d > 0$ such that the last five inequalities hold for all 
$\| g \|_\b \leq \d$, namely we fix $\d = 1 / \max \{ K_1$, $K_2$, $C_* K_1 (1+A) / \d_1$, 
$C_* M_{123}(a_2 - \mu) L_{456}(a_2) (K_1 + K_3) \}$.
This completes the proof of \eqref{Olli.9}, \eqref{Olli.10}, \eqref{Olli.11}, \eqref{reve}.

\bigskip

\noindent
\textbf{Convergence of the scheme.}
The same argument used in \eqref{3005.1}, \eqref{3005.2}, \eqref{3005.3}, \eqref{0106.1} 
proves that $(u_n)$ is a Cauchy sequence in $E_\a$.  
Hence $u_n$ converges to a limit $u \in E_\a$, with $\| u \|_\a \leq K_4 \| g \|_\b$.

We prove the convergence of the scheme. 
By \eqref{new.1} and \eqref{new.24} one proves by induction that
\[ 
\sum_{j=0}^k (e_j + y_j) = e_k + r_k, \quad \text{where} \ \ 
r_k := (I - S_k) \sum_{j=0}^{k-1} e_j,
\quad \forall k \geq 1.	
\] 
Hence, by \eqref{Olli.8} and \eqref{new.2}, 
recalling that $\Phi'(v_j) \Psi(v_j)$ is the identity map,
one has
\[ 
\Phi(u_{k+1}) - \Phi(u_0)
= \sum_{j=0}^{k}[\Phi(u_{j+1}) - \Phi(u_j)] 
= \sum_{j=0}^{k} (e_j + g_j + y_j) 
= G_k + e_k + r_k
\]
where $G_k := \sum_{j=0}^{k} g_j = S_{k+1} g$. 
By \eqref{S3}, \eqref{S1}, $\| G_k - g \|_b \to 0$ as $k \to \infty$, for all $b \in [0,\b)$.
By \eqref{new.9}, \eqref{Olli.9}, \eqref{Olli.11} and \eqref{reve},
$\| e_j \|_{\a - \mu} \leq M \, 2^{j(a_1 - \a)}$ for some $M>0$, 
and the series $\sum_{j=0}^\infty \| e_j \|_{\a-\mu}$ converges.  
By \eqref{S3}, for all $\rho \in [0,\a-\mu)$ we have 
\begin{equation}
\| r_k \|_\rho 
\leq \sum_{j=0}^{k-1} \| (I - S_k) e_j	\|_\rho
\leq \sum_{j=0}^{k-1} C_\rho 2^{-k(\a - \mu - \rho)} \| e_j	\|_{\a-\mu}
\leq C_\rho M 2^{-k(\a - \mu - \rho)},
\end{equation}
so that $\| r_k \|_\rho \to 0$ as $k \to \infty$.
We have proved that $\| \Phi(u_k) - \Phi(u_0) - g \|_\rho \to 0$ as $k \to \infty$ 
for all $\rho$ in the interval $0 \leq \rho < \min \{ \a - \mu, \b \}$. 
Since $u_k \to u$ in $E_\a$, 
it follows that $\Phi(u_k) \to \Phi(u)$ in $F_{\a-\mu}$. 
This completes the proof of the first part of the theorem.

\bigskip

\noindent
\textbf{Higher regularity.} 
It remains to prove the last part of the theorem. 
Let $c > 0$. 
Assume that \eqref{Phi sec} holds for all $a \in [0, a_2 + c - \mu]$,
and that \eqref{tame in NM} holds for all $a \in [a_1, a_2 + c]$.
Assume that $g \in F_{\b+c}$, with \eqref{0406.1}.
By \eqref{2705.3}, 
\begin{equation} \label{0606.3}
\| g_j \|_b \leq C_{b,c} \, 2^{j(b-\b-c)} \| g_j \|_{\b+c} \quad \forall b \geq 0	
\end{equation}
(namely \eqref{2705.2} holds for $b \in [0,\infty)$, with $\b$ replaced by $\b+c$).

\noindent
$\bullet$ \emph{Improved estimates.}
Using \eqref{S2}, \eqref{new.16}, \eqref{2705.3}, \eqref{new.18}, 
and \eqref{3110.1}, we have 
\begin{align} 
\| y_{k+1} \|_b 
& \leq C_b K_1 (K_1 + K_3) M_{123}(a_2 - \mu) \| g \|_\b^2 \, 2^{(k+1)(b-\b-\g)} 
\notag \\ 
& \leq C_b \| g \|_\b \, 2^{(k+1)(b-\b-\g)} \quad \forall b \geq 0
\label{1810.1}
\end{align}
(namely \eqref{new.8} holds for $b \in [0, \infty)$, with $C$ replaced by $C_b$, 
then we use \eqref{3110.1}, recalling that $K_1 = C_* L_{456}(a_2)$). 
Using \eqref{S2}, \eqref{Olli.10} and \eqref{3110.1}, we have
\begin{equation} \label{1810.2}
\| v_j \|_a 
\leq C_a K_2 \| g \|_\b \, 2^{j(a-\a)} 
\leq C_a 2^{j(a-\a)} 
\quad \forall a \geq a_1 + \b
\end{equation}
(namely \eqref{Olli.10} holds for $a \in [a_1 + \b, \infty)$, 
with $K_2$ replaced by $C_a K_2$, then use \eqref{3110.1}).
By \eqref{Olli.8}, \eqref{tame in NM} (which now holds for $a \in [a_1, a_2+c]$), 
\eqref{0606.3}, \eqref{1810.1}, \eqref{1810.2}, 
and \eqref{2705.2} for the term containing $L_6(a) \| g_k \|_0$,
we deduce that, for all $k \geq 0$, 
\begin{align} 
\| h_k \|_a 
& \leq L_{45}(a) \big( C_{a,c} \| g_k \|_{\b+c} \, 2^{k(a - \a - c)}
+ C_a \| g \|_\b \, 2^{k(a - \a - \g)} \big) + L_6(a) \mC 2^{-k\b} \xi_k 
\notag \\ 
& \leq L_{45}(a) C_{a,c} 2^{k(a - \a - \lm)} \eta_k + L_6(a) \mC 2^{-k\b} \psi_k
\quad \forall a \in [a_1, a_2+c],
\label{0606.1}
\end{align}
where $L_{45} := L_4 + L_5$, 
$\mC$ is the sum of the two constants $C_b$ at $b=0$ 
appearing in \eqref{2705.2} and \eqref{1810.1},
$\xi_k$ has been defined above as $\xi_k = \| g \|_\b \, 2^{- k\g} + \| g_k \|_\b$,
\begin{equation} \label{0111.7}
\eta_k := \| g_k \|_{\b+c} + \| g \|_{\b+c} \, 2^{-k \g / 2}, \quad 
\psi_k := \| g_k \|_\b + \| g \|_\b \, 2^{-k \g / 2}, \quad 
\lm := \frac{c}{N}\,,
\end{equation}
and $N$ is the smallest positive integer that is $\geq 2c/\g$
(so that $\lm \leq \min \{ c, \g/2 \}$ and $N\lm = c$).
For $a = a_1$, by \eqref{Olli.8}, \eqref{tame in NM}, 
\eqref{0606.3} (which here we use also for the term containing $L_6(a_1) \| g_k \|_0$), 
\eqref{1810.1} and \eqref{1810.2}, since $- \b < a_1 - \a$, we obtain
\begin{equation} \label{2610.2}
\| h_k \|_{a_1} \leq C_c L_{456}(a_1) 2^{k(a_1 - \a - \lm)} \eta_k.
\end{equation}

\noindent
$\bullet$ \emph{Finite induction.} 
If $N=1$, then \eqref{0606.1} gives \eqref{0211.3} below.
If, instead, $N \geq 2$, we repeat the argument and prove recursively for $n=1,\ldots,N$ the following bounds: 
for all $k \geq 0$, all $a \in [a_1, a_2 + c]$,
\begin{align} \label{2710.1} 
\| h_k \|_a 
& \leq 2^{k (a -\a - n\lm)} ( \mA_n(a) \psi_k + \mB_n(a) \eta_k) 
+ 2^{-k \b} L_6(a) \mC \psi_k,
\\ 
\| h_k \|_{a_1} 
& \leq 2^{k (a_1 -\a - n\lm)} (\mE_n \psi_k + \mF_n \eta_k),
\label{2710.2}
\end{align}
where the coefficients $\mA_n(a), \mB_n(a), \mE_n, \mF_n$ are defined recursively,
and $\mC$ has been defined above as the sum of the two constants $C_b$ at $b=0$ appearing in \eqref{2705.2} and \eqref{1810.1}. 
Estimates \eqref{0606.1} and \eqref{2610.2} give \eqref{2710.1}, \eqref{2710.2} for $n=1$ with 
\begin{equation} \label{0211.1}
\mA_1(a) = \mE_1 = 0, \quad 
\mB_1(a) = L_{45}(a) C_{a,c}, \quad 
\mF_1 = L_{456}(a_1) C_c.
\end{equation}
Suppose that \eqref{2710.1}-\eqref{2710.2} hold for some $n \in [1,N-1]$. 
We have to prove that they also hold for $n+1$.
By \eqref{2710.1}, since $\psi_k \leq C \| g \|_\b$, 
$\eta_k \leq C_c \| g \|_{\b+c}$,  and $a_2 + c - \a - n\lm > 0$, 
\begin{align} \label{2510.1}
\| u_k \|_{a_2+c} 
& \leq \sum_{j=0}^{k-1} \| h_j \|_{a_2+c} 
\notag \\ 
& \leq 2^{k (a_2 + c -\a - n\lm)} ( \tilde \mA_n C \| g \|_\b + \tilde \mB_n C_c \| g \|_{\b+c} )
+ \tilde L_6 \mC C \| g \|_\b,
\end{align}
where $\tilde \mA_n := \mA_n(a_2+c)$, $\tilde \mB_n := \mB_n(a_2 + c)$, 
$\tilde L_6 := L_6(a_2+c)$. 
By \eqref{S1}, $\| v_k \|_{a_2+c} \leq C_c \| u_k \|_{a_2+c}$. 
Therefore $v_k$ satisfies the same bound \eqref{2510.1} as $u_k$,
and, by triangle inequality, $\| v_k - u_k \|_{a_2 + c}$ also does. 

By assumption, \eqref{Phi sec} holds for $a \in [0, a_2 + c - \mu]$. 
Therefore \eqref{new.9} also holds for $a$ in the same interval, 
and it can be used to estimate $\| e_j \|_{a_2 + c - \mu}$.
Using \eqref{Olli.9}, \eqref{Olli.11}, 
\eqref{3110.1} for the ``low norm'' factors $\| h_j \|_{a_1}$, $\| v_j - u_j \|_{a_1}$, 
and \eqref{2710.1}, \eqref{2510.1} for the ``high norm'' factors 
$\| h_j \|_{a_2 + c}$, $\| u_j \|_{a_2 + c}$, $\| v_j \|_{a_2 + c}$, $\| v_j-u_j \|_{a_2 + c}$, 
we obtain
\begin{align} 
\| e_j \|_{a_2 + c - \mu} 
& \leq 2^{j(a_1 + a_2 - 2\a + c - n\lm)} 
\big\{ \tilde \mA_n \tilde M_{12} C \| g \|_\b + \tilde \mB_n \tilde M_{12} C_c \| g \|_{\b+c} \big\} 
\notag \\ 
& \quad \ 
+ 2^{j(a_1 - \a)} 
\big\{ \tilde L_6 \mC \tilde M_{12} C \| g \|_\b 
+ \tilde M_3 K_1 \| g \|_\b \big\}
\label{1810.5}
\end{align}
where $\tilde M_i := M_i (a_2 + c - \mu)$, $i = 1,2,3$, and $\tilde M_{12} := \tilde M_1 + \tilde M_2$. 

By \eqref{new.9}, \eqref{Olli.9}, \eqref{Olli.11},  
\eqref{3110.1} we have $\| e_j \|_0 \leq 2^{j(a_1 -\a)} \| h_j \|_{a_1} M_{123}(0)$. 
Hence, by \eqref{2710.2}, 
\begin{equation} \label{1810.4}
\| e_j \|_0 \leq 2^{j (2a_1 - 2\a - n\lm)} \big\{ \mE_n M_{123}(0) \psi_j 
+ \mF_n M_{123}(0) \eta_j \big\}.
\end{equation} 

By \eqref{S2}, $\| S_{k+1} e_k \|_b \leq C_b 2^{(k+1)b} \| e_k \|_0$ for all $b \geq 0$, 
and therefore, using \eqref{1810.4}, 
we obtain an estimate for $\| S_{k+1} e_k \|_b$ for all $b \geq 0$.
By \eqref{2705.3}, for all $b \geq 0$,
\[
\sum_{j=0}^{k-1} \| R_k e_j \|_b 
\leq C_{b,c} 2^{k(b - a_2 - c + \mu)} \sum_{j=0}^{k-1} \| e_j \|_{a_2 + c - \mu},
\]
and therefore, using \eqref{1810.5} and the fact that $(a_1 + a_2 - 2\a + c - n\lm) > 0$, 
we get an estimate for $\| R_k \sum_{j=0}^{k-1} e_j \|_b$ for all $b \geq 0$.
Recalling \eqref{new.1}, we deduce that, for all $k \geq 0$, 
\begin{align} 
\| y_{k+1} \|_b 
& \leq 
2^{(k+1)(b - a_2 - c + \mu)} \big\{ \tilde L_6 \mC \tilde M_{12} C_{b,c} \| g \|_\b 
+ \tilde M_3 C_{b,c} K_1 \| g \|_\b \big\}
\notag \\ & \quad \ 
+ 2^{(k+1)(b + 2a_1 - 2\a - n\lm)} \big\{ \mE_n M_{123}(0) C_b \psi_k 
+ \mF_n M_{123}(0) C_b \eta_k 
\notag \\ & \quad \ 
+ \tilde \mA_n \tilde M_{12} C_{b,c} \| g \|_\b
+ \tilde \mB_n \tilde M_{12} C_{b,c} \| g \|_{\b+c} \big\}
\qquad \forall b \geq 0.
\label{0606.2}
\end{align}
The exponents in \eqref{0606.2} satisfy 
$(b - a_2 - c + \mu) \leq (b + 2a_1 - 2\a - n\lm)$, 
because $a_1 + a_2 - 2 \a > 0$ and $c = N \lm > n \lm$.  
Moreover, $(b + 2a_1 - 2\a - n\lm) \leq (b - \b - (n+1) \lm - (\g/2))$ 
because $\lm \leq \g/2$ and $2a_1 - 2 \a + \b + \g \leq 0$.  
Hence, for all $k \geq 0$,  
\begin{equation} \label{0111.4}
\| y_k \|_b \leq 2^{k (b - \b - (n+1) \lm - (\g/2))} \, C_{b,c} Y_n 
\quad \forall b \geq 0,
\end{equation}
where
\begin{align}
Y_n & := \big\{ \tilde \mA_n \tilde M_{12} + \tilde L_6 \mC \tilde M_{12}
+ K_1 \tilde M_3
+ \mE_n M_{123}(0) \big\} \| g \|_\b
\notag \\ & \qquad 
+ \big\{ \tilde \mB_n \tilde M_{12} + \mF_n M_{123}(0) \big\} \| g \|_{\b+c}.
\label{0111.3}
\end{align}

By \eqref{Olli.8} and \eqref{tame in NM} we estimate $\| h_k \|_a$ for $a \in [a_1, a_2 + c]$. 
Since $c = N \lm \geq (n+1)\lm$, 
using \eqref{0606.3}, \eqref{1810.2} for $L_4(a) \| g_k \|_{a+\b-\a} + L_5(a) \| v_k \|_{a+\b} \| g_k \|_0$, 
and \eqref{2705.2} for $L_6(a) \| g_k \|_0$, we get, for all $a \in [a_1, a_2 + c]$,
\begin{equation} \label{0111.1}
\| \Psi(v_k) g_k \|_a \leq 2^{k(a-\a-(n+1)\lm)} L_{45}(a) C_{a,c} \| g_k \|_{\b+c} 
+ 2^{-k\b} L_6(a) C \| g_k \|_\b .
\end{equation}
Using \eqref{0111.4}, \eqref{1810.2} for $L_4(a) \| y_k \|_{a+\b-\a} + L_5(a) \| v_k \|_{a+\b} \| y_k \|_0$	
and \eqref{1810.1} for $L_6(a) \| y_k \|_0$, we get, for all $a \in [a_1, a_2 + c]$,
\begin{equation} \label{0111.5}
\| \Psi(v_k) y_k \|_a \leq 
2^{k(a-\a-(n+1)\lm)} L_{45}(a) C_{a,c} Y_n 2^{-k \g/2} + 2^{-k\b} L_6(a) C \| g \|_\b 2^{-k\g}.
\end{equation}
Recalling that $K_1 = C_* L_{456}(a_2)$ and the definition \eqref{0111.7} of $\psi_k, \eta_k$, 
the sum of \eqref{0111.1} and \eqref{0111.5} gives \eqref{2710.1} at $n+1$, with
\begin{align} \label{0111.8}
\mA_{n+1}(a) & = L_{45}(a) C_{a,c} (\tilde \mA_n \tilde M_{12} + \tilde L_6 \tilde M_{12}
+ L_{456}(a_2) \tilde M_3 + \mE_n M_{123}(0)),
\\
\label{0111.9}
\mB_{n+1}(a) & = L_{45}(a) C_{a,c} (1 + \tilde \mB_n \tilde M_{12} + \mF_n M_{123}(0)).
\end{align}
Using \eqref{1810.2}, \eqref{0606.3} also for the term $L_6(a_1) \| g_k \|_0$, we get 
\begin{equation} \label{0111.21}
\| \Psi(v_k) g_k \|_{a_1} \leq 2^{k(a_1-\a-(n+1)\lm)} L_{456}(a_1) C_c \| g_k \|_{\b+c}. 
\end{equation}
Using \eqref{0111.4}, \eqref{1810.2} also for the term $L_6(a_1) \| y_k \|_0$, we get
\begin{equation} \label{0111.25}
\| \Psi(v_k) y_k \|_{a_1} \leq 2^{k(a_1 -\a-(n+1)\lm)} L_{456}(a_1) C_c Y_n 2^{-k \g/2}.
\end{equation}
The sum of the last two bounds gives \eqref{2710.2} at $n+1$, with
\begin{align} \label{0111.11}
\mE_{n+1} & = L_{456}(a_1) C_c (\tilde \mA_n \tilde M_{12} + \tilde L_6 \tilde M_{12} 
+ L_{456}(a_2) \tilde M_3 + \mE_n M_{123}(0)) 
\\
\label{0111.12}
\mF_{n+1} & = L_{456}(a_1) C_c (1 + \tilde \mB_n \tilde M_{12} + \mF_n M_{123}(0)).
\end{align}
Let 
\begin{equation} \label{0111.32}
Z := 	L_{456}(a_1) C_c M_{123}(0) + \tilde L_{45} \tilde C_c \tilde M_{12}, \quad 
X := \tilde L_6 \tilde M_{12} + L_{456}(a_2) \tilde M_3,
\end{equation}
where the constant $C_c$ in \eqref{0111.32} is the one of \eqref{0111.11}-\eqref{0111.12}, 
and the constant $\tilde C_c$ is the constant $C_{a,c}$ of \eqref{0111.8}-\eqref{0111.9}
evaluated at $a = a_2 + c$. 
By induction, the recursive system \eqref{0111.8}, \eqref{0111.9}, \eqref{0111.11}, \eqref{0111.12} with the initial values \eqref{0211.1} gives
\begin{alignat}{2} \label{0111.30}
\mA_n(a) & = L_{45}(a) C_{a,c} X \sum_{j=0}^{n-2} Z^j, 
\quad & \quad 
\mB_n(a) & = L_{45}(a) C_{a,c} \sum_{j=0}^{n-1} Z^j, 
\\ 
\mE_n & = L_{456}(a_1) C_c X \sum_{j=0}^{n-2} Z^j, 
\quad & \quad 
\mF_n & = L_{456}(a_1) C_c \sum_{j=0}^{n-1} Z^j
\label{0111.31}
\end{alignat}
for all $n \geq 2$. 
The iteration ends at $n = N$, and, since $N \lm = c$, we obtain for all $k \geq 0$
\begin{equation} \label{0211.3} 
\| h_k \|_a \leq 2^{k (a -\a - c)} (\mA_N(a) \psi_k + \mB_N(a) \eta_k) 
+ 2^{-k \b} L_6(a) \mC \psi_k 
\quad \forall a \in [a_1, a_2+c].
\end{equation}

\noindent
$\bullet$ \emph{Convergence in high norm.} 
The argument used in \eqref{3005.1}-\eqref{0106.1} 
(now with $a_1 + c, \a + c, a_2 + c$ instead of $a_1, \a, a_2$,
and bound \eqref{0211.3} instead of \eqref{Olli.9}) 
proves that $(u_n)$ is a Cauchy sequence in $E_{\a+c}$, 
and its limit $u$ satisfies 
\begin{equation} \label{0211.4}
\| u \|_{\a+c} \leq C(c) \{ (\tilde L_6 + \tilde \mA_N) (1 + A) \| g \|_\b 
+ \tilde \mB_N (1 + A_c) \| g \|_{\b+c} \}
\end{equation}
for some constant $C(c)$ depending on $c$. 
The proof of Theorem \ref{thm:1} is complete.
\qed

\bigskip

\begin{remark} \label{rem:11.09.2017}
In \cite{Olli}, the bound corresponding to \eqref{Olli.9} (estimate (9) in \cite{Olli})
is $\| \dot{u}_j \|_a \leq C_1 \| g \|_{F_\b'} \theta_j^{a - \a - 1}$ for all $a \in [a_1, a_2]$,
where $F_\b'$ is the weak space whose definition is recalled in Remark \ref{rem:2707.2}.
In our notation with $\theta_j = 2^j$ this corresponds to 
$h_j = 2^j \dot{u}_j$ and
\begin{equation} \label{sept.9}
\| h_j \|_a \leq C_1 \| g \|_{F_\b'} 2^{j(a - \a)} 
\quad \forall a \in [a_1, a_2].
\end{equation}
Also, in \cite{Olli} the bound corresponding to \eqref{reve} (estimate (12) in \cite{Olli}) 
is 
\begin{equation} \label{sept.12}
\| u_j \|_{E_\a'} \leq C' C_1 \| g \|_{F_\b'}. 
\end{equation}
Estimate \eqref{sept.9} at the regularity threshold $a = \a$ only implies \eqref{sept.12},
and therefore \eqref{sept.9} is sufficient to deduce that the solution 
$u = \sum_{j=0}^\infty h_j$ belongs to the weak space $E_\a'$, 
but it is not sufficient to prove that $u \in E_\a$. 
For this reason, when the datum $g \in F_\b$, 
the implicit function theorem in \cite{Olli} and the one in \cite{Geodesy}
give a solution $u$ of the equation $\Phi(u) = \Phi(0) + g$ that only belongs 
to the weak space $E_\a'$, which, in the Sobolev case, is larger than $E_\a$.

The solution $u$ given by Theorem \ref{thm:1}, instead, belongs to $E_\a$ 
when the datum $g \in F_\b$ satisfies the ``orthogonality assumption'' \eqref{2705.1}.
To obtain this sharp regularity we use a stronger version of \eqref{sept.9}-\eqref{sept.12}
given by \eqref{Olli.9} and \eqref{reve}. 
Note that the factor $(\| g \|_\b 2^{-j\g} + \| g_j \|_\b)$ in \eqref{Olli.9} 
(see also $\xi_k$ in \eqref{3005.1} and $\eta_k, \psi_k$ in \eqref{0111.7})
has a stronger summability property than the corresponding factor $\| g \|_{F_\b'}$ of \eqref{sept.9}
 --- at the threshold $a = \a$ the right hand side of \eqref{Olli.9} is a sequence in $\ell^2$, 
while the right hand side of \eqref{sept.9} is only in $\ell^\infty$.

However, it is not trivial to deduce \eqref{reve} from \eqref{Olli.9} 
(remember that $h_j$ is not the $j$-th dyadic block of $u$).
This is the point where we apply a trick inspired by paradifferential calculus
(see for example the proof of Proposition 4.1.13 on page 53 of M\'etivier \cite{Metivier}).  
To estimate $\| u_{k+1} \|_\a$, we first use the dyadic decomposition
$u_{k+1} = \sum_{j=0}^\infty R_j u_{k+1}$. Then we use the identity $u_{k+1} = \sum_{n=0}^k h_n$ 
(see the recursive scheme \eqref{Olli.8}), and estimate the norm $\| R_j h_n \|_\a$ 
of each dyadic block of each component. 
The estimate is performed according to the frequency localization: 
the terms $R_j h_n$ with $n \leq \min \{k,j\}$ 
(where the iteration index $n$ is smaller than the frequency localization $j$)
are collected in the sum $\e_j''$ in \eqref{3005.3} and are estimated 
using the high norm $a_2$, 
while possible terms with $n > j$ 
(where the iteration index is larger than the frequency localization)
are collected in the sum $\e_j'$ in \eqref{3005.3} 
and are estimated using the low norm $a_1$. 
Then the dyadic decomposition, and the fact that $(\xi_n) \in \ell^2$, 
are used to estimate the $\ell^2$ norm of the corresponding sequence (see \eqref{0106.1}).
Finally the orthogonality assumption \eqref{2705.4} for the dyadic decomposition $(R_j)$ 
gives \eqref{reve}.
\end{remark}

\section{Application to quasi-linear perturbations of KdV}
\label{sec:appl}

We use Theorem \ref{thm:1} to improve the regularity in the results of exact controllability and local well-posedness for the Cauchy problem of quasi-linear perturbations of KdV obtained in \cite{BFH}.

We consider  equations of the form
\begin{equation} \label{i1}
u_t + u_{xxx} + \ttf(x,u, u_x, u_{xx}, u_{xxx}) = 0
\end{equation}
where the nonlinearity $\ttf(x,u,u_x, u_{xx}, u_{xxx})$ 
is at least quadratic around $u=0$, namely the real-valued function 
$\ttf : \T \times \R^4 \to \R$ satisfies 
\begin{equation} \label{i2}
|\ttf(x, z_0, z_1, z_2, z_3)| \leq C |z|^2 
\quad \forall z = (z_0, z_1, z_2, z_3) \in \R^4, \ |z| \leq 1.
\end{equation}
We assume that the dependence of $\ttf$ on $u_{xx}, u_{xxx}$ is Hamiltonian, 
while no structure is required on its dependence on $u, u_x$. More precisely, we assume that 
\begin{equation} \label{i2.1}
\ttf(x,u,u_x, u_{xx}, u_{xxx}) 
= \ttf_1(x,u,u_x, u_{xx}, u_{xxx}) + \ttf_0(x,u,u_x)
\end{equation}
where 
\begin{equation}  \label{i6}
\begin{aligned} 
& \ttf_1(x,u,u_x, u_{xx}, u_{xxx}) 
= \pa_x \{ (\pa_u \mF)(x, u, u_x) \} 
- \pa_{xx} \{ (\pa_{u_x} \mF)(x, u, u_x) \}
\\
& \text{for some function $\mF : \T \times \R^2 \to \R$.}
\end{aligned}
\end{equation}
Note that the case $\ttf = \ttf_1$, $\ttf_0 = 0$ corresponds to the Hamiltonian equation 
$\pa_t u = \pa_x \gr H(u)$ where the Hamiltonian is 
\begin{equation} \label{i3}
H(u) = \frac12 \int_\T u_x^2 \, dx + \int_\T \mF(x,u,u_x) \, dx
\end{equation}
and $\gr$ denotes the $L^2(\T)$-gradient. The unperturbed KdV is the case $\mF = - \frac16 u^3$.

\begin{theorem}[Exact controllability] \label{thm:control}
Let $T>0$, and let $\om \subset \T$ be a nonempty open set. 
There exist positive universal constants $r_1,s_1$ such that, 
if $\mN$ in \eqref{i1} is of class $C^{r_1}$ in its arguments 
and satisfies \eqref{i2}, \eqref{i2.1}, \eqref{i6}, 
then there exists a positive constant $\d_*$ depending on $T,\om,\mN$
with the following property. 

Let $u_{in}, u_{end} \in H^{s_1}(\T,\R)$ 
with 
\[
\| u_{in} \|_{s_1} + \| u_{end} \|_{s_1} \leq \d_*.
\] 
Then there exists a function $f(t,x)$ satisfying
\[
f(t,x) = 0 \quad \text{for all $x \notin \om$, for all $t \in [0,T]$,}
\]
belonging to $C([0,T],H^{s_1}_x) 
\cap C^1([0,T],H^{s_1-3}_x) 
\cap C^2([0,T],H^{s_1-6}_x)$
such that the Cauchy problem
\begin{equation} \label{i9}
\begin{cases}
u_t + u_{xxx} + \ttf(x,u,u_x, u_{xx}, u_{xxx}) = f 
\quad \forall (t,x) \in [0,T] \times \T \\
u(0,x) = u_{in}(x) 
\end{cases}
\end{equation}
has a unique solution $u(t,x)$ belonging to $C([0,T], H^{s_1}_x) \cap C^1([0,T], H^{s_1-3}_x)
\cap C^2([0,T],H^{s_1-6}_x)$,
which satisfies 
\begin{equation} \label{i10}
u(T,x) = u_{end}(x),
\end{equation}
and
\begin{multline} \label{stimetta}
\| u,f \|_{C([0,T],H^{s_1}_x)} + \| \pa_t u, \pa_t f \|_{C([0,T],H^{s_1-3}_x)} 
+ \| \pa_{tt} u, \pa_{tt} f \|_{C([0,T],H^{s_1-6}_x)} 
\\
\leq C_{s_1} (\| u_{in} \|_{s_1} + \| u_{end} \|_{s_1})
\end{multline}
for some $C_{s_1} > 0$ depending on $s_1,T,\om,\mN$.

Moreover, the universal constant $\tau_1:= r_1 - s_1 > 0$ has the following property.
For all $r \geq r_1$, all $s \in [s_1, r - \tau_1]$, if, in addition to the previous assumptions, $\mN$ is of class $C^r$ and $u_{in}, u_{end} \in H^s_x$, then $u,f$
belong to $C([0,T], H^s_x) \cap C^1([0,T], H^{s-3}_x)
\cap C^2([0,T],H^{s-6}_x)$ and \eqref{stimetta} holds with $s$ instead of $s_1$.
\end{theorem}

\begin{theorem}[Local existence and uniqueness] \label{thm:byproduct}
There exist positive universal constants $r_0,s_0$ such that, 
if $\mN$ in \eqref{i1} is of class $C^{r_0}$ in its arguments 
and satisfies \eqref{i2}, \eqref{i2.1}, \eqref{i6}, 
then the following property holds.
For all $T > 0$ there exists $\d_* > 0$ such that
for all $u_{in} \in H^{s_0}_x$ satisfying
\begin{equation} \label{i12}
\| u_{in} \|_{s_0} \leq \d_* \,,
\end{equation}
the Cauchy problem 
\begin{equation} \label{i11}
\begin{cases}
u_t + u_{xxx} + \ttf (x,u,u_x, u_{xx}, u_{xxx}) = 0, 
\qquad (t,x) \in [0,T] \times \T \\
u(0,x) = u_{in}(x) 
\end{cases}
\end{equation}
has one and only one solution $u \in C([0,T], H^{s_0}_x) \cap C^1([0,T], H^{s_0-3}_x) \cap C^2([0,T], H^{s_0-6}_x)$.
Moreover
\begin{equation} \label{stimetta bis}
\| u \|_{C([0,T],H^{s_0}_x)} + \| \pa_t u \|_{C([0,T],H^{s_0-3}_x)} 
+ \| \pa_{tt} u \|_{C([0,T],H^{s_0-6}_x)} 
\leq C_{s_0} \| u_{in} \|_{s_0}
\end{equation}
for some $C_{s_0} > 0$ depending on $s_0,T,\mN$.

Moreover the universal constant $\tau_0 := r_0 - s_0 > 0$ has the following property. For all $r \geq r_0$, all $s \in [s_0, r - \tau_0]$, if, in addition to the previous assumptions, $\mN$ is of class $C^r$ and $u_{in} \in H^s_x$, then $u$ belongs to $C([0,T], H^s_x) \cap C^1([0,T], H^{s-3}_x)
\cap C^2([0,T],H^{s-6}_x)$ and \eqref{stimetta bis} holds with $s$ instead of $s_0$.
\end{theorem} 

\noindent
{\bf Proof of Theorem \ref{thm:control}.} 
Define
\begin{equation} \label{ge2}
P(u) := u_t + u_{xxx} + \ttf(x,u,u_x,u_{xx}, u_{xxx}).
\end{equation}
and
\begin{equation} \label{ge3}
\Phi(u,f) := \begin{pmatrix} 
P(u) - \chi_\om f \\
u(0) \\
u(T) \end{pmatrix}
\end{equation}
so that the problem
\begin{equation} \label{2507.1}
\begin{cases}
u_t + u_{xxx} + \ttf(x,u,u_x, u_{xx}, u_{xxx}) = f 
\quad \forall (t,x) \in [0,T] \times \T \\
u(0,x) = u_{in}(x) \\
u(T,x) = u_{end}(x)
\end{cases}
\end{equation}
is written as $\Phi(u,f) = (0, u_{in}, u_{end})$.
The linearized operator $\Phi'(u,f)[h,\ph]$ at the point $(u,f)$ in the direction $(h,\ph)$ is
\begin{equation} \label{ge4}
\Phi'(u,f)[h,\ph] := \begin{pmatrix} 
P'(u)[h] - \chi_\om \ph \\
h(0) \\
h(T) \end{pmatrix}.
\end{equation}

We define the scales of Banach spaces 
\begin{equation} \label{def Es}
E_s := X_s \times X_s, \quad 
X_s := C([0,T], H^{s+6}_x) \cap C^1([0,T], H^{s+3}_x) \cap C^2([0,T], H^s_x)
\end{equation}
and 
\begin{equation} \label{def Fs}
F_s := \{ g = (g_1, g_2, g_3) : g_1 \in C([0,T], H^{s+6}_x) \cap C^1([0,T], H^s_x), 
g_2, g_3 \in H^{s+6}_x \}
\end{equation}
equipped with the norms 
\begin{equation} \label{def norm Es}
\| u,f \|_{E_s} := \| u \|_{X_s} + \| f \|_{X_s}, \quad 
\| u \|_{X_s} := \| u \|_{T,s+6} + \| \pa_t u \|_{T,s+3} + \| \pa_{tt} u \|_{T,s}
\end{equation}
and
\begin{equation} \label{def norm Fs}
\| g \|_{F_s} := \| g_1 \|_{T,s+6} + \| \pa_t g_1 \|_{T,s} + \| g_2, g_3 \|_{s+6}.
\end{equation}
In Theorem 4.5 of \cite{BFH}, the following right inversion result for the linearized operator in \eqref{ge4} is proved.

\begin{proposition} \label{thm:inv}
Let $T>0$, and let $\om \subset \T$ be an open set.
There exist two universal constants $\t,\s \geq 3$ 
and a positive constant $\d_*$ depending on $T,\om$ with the following property.

Let $s \in [0, r-\t]$, 
where $r$ is the regularity of the nonlinearity $\mN$. 
Let $g = (g_1, g_2, g_3) \in F_s$, 
and let $(u,f) \in E_{s+\s}$, 
with $\| u \|_{X_\s} \leq \d_*$. 
Then there exists $(h,\ph) := \Psi(u,f)[g] \in E_s$ such that 
\begin{equation} \label{3009}
P'(u)[h] - \chi_\om \ph = g_1, \quad 
h(0) = g_2, \quad 
h(T) = g_3,
\end{equation}
and 
\begin{equation}  \label{2809}
\| h,\ph \|_{E_s} \leq C_s \big( \| g \|_{F_s} + \| u \|_{X_{s+\s}} \| g \|_{F_0} \big)
\end{equation}
where $C_s$ depends on $s,T,\om$.
\end{proposition}

We define the smoothing operators $S_j$, $j = 0,1,2,\ldots$ as  
\[ 
S_j u (x) := \sum_{|k| \leq 2^j} \widehat u_k \, e^{i k x}
\qquad \text{where} \quad 
u(x) = \sum_{k \in \Z} \widehat u_k \, e^{i k x} 
\] 
The definition of $S_j$ extends in the obvious way to functions 
$u(t,x) = \sum_{k \in \Z} \widehat u_k(t) \, e^{i kx}$ 
depending on time. 
Since $S_j$ and $\pa_t$ commute, the smoothing operators $S_j$ are defined 
on the spaces $E_s$, $F_s$ defined in \eqref{def Es}-\eqref{def Fs} by setting 
$S_j(u,f) := (S_j u, S_j f)$ and similarly on $g = (g_1, g_2, g_3)$. 
One easily verifies that $S_j$ satisfies \eqref{S0}-\eqref{S4} and \eqref{2705.4} 
on $E_s$ and $F_s$. 

By \eqref{ge3}, observe that $\Phi(u,f) := (P(u) - \chi_\om f, \, u(0), \, u(T) )$ 
belongs to $F_s$ 
when $(u,f) \in E_{s+3}$, $s \in [0, r-6]$, 
with $\| u \|_{T,4} \leq 1$. 
Its second derivative in the directions $(h,\ph)$ and $(w,\psi)$ is
\[
\Phi''(u,f)[(h, \ph), (w,\psi)] 
= \begin{pmatrix} P''(u)[h, w] \\ 0 \\ 0 \end{pmatrix}.
\]
For $u$ in a fixed ball $\| u \|_{X_1} \leq \d_0$, with $\d_0$ small enough, 
we estimate
\begin{equation} \label{stima Phi''}
\| P''(u)[h,w] \|_{F_s} 
\lesssim_s  \big( \| h \|_{X_1} \| w \|_{X_{s+3}} 
+ \| h \|_{X_{s+3}} \| w \|_{X_1}
+ \| u \|_{X_{s+3}} \| h \|_{X_1} \| w \|_{X_1} \big)
\end{equation}
for all $s \in [0, r-6]$. 
We fix $V = \{ (u,f) \in E_3 : \| (u,f) \|_{E_3} \leq \d_0 \}$,
$\d_1 = \d_*$,  
\begin{equation} \label{param.1}
a_0 = 1, \quad
\mu = 3, \quad
a_1 = \s, \quad
\a = \b > 2 \s, \quad
a_2 > 2 \a - a_1 \quad
\end{equation}
where $\d_*, \s, \t$ are given by Proposition \ref{thm:inv}, 
and $r \geq r_1 := a_2 + \t$ is the regularity of $\mN$. 
The right inverse $\Psi$ in Proposition \ref{thm:inv} satisfies the assumptions of Theorem \ref{thm:1}.
Let $u_{in}, u_{end} \in H^{\b+6}_x$, with $\| u_{in}, u_{end} \|_{H^{\b+6}_x}$ small enough.
Let $g := (0, u_{in}, u_{end})$, so that $g \in F_\b$ and $\| g \|_{F_\b} \leq \d$.
Since $g$ does not depend on time, it satisfies \eqref{2705.1}.

Thus by Theorem \ref{thm:1} there exists a solution $(u,f) \in E_\a$ of the equation
$\Phi(u,f) = g$, with $\| u,f \|_{E_\a} \leq C \| g \|_{F_\b}$ 
(and recall that $\b = \a$). 
We fix $s_1 := \a + 6$, and \eqref{stimetta} is proved. 

We have found a solution $(u,f)$ of the control problem \eqref{2507.1}. 
Now we prove that $u$ is the unique solution of the Cauchy problem \eqref{i9}, 
with that given $f$. 
Let $u,v$ be two solutions of \eqref{i9} in $E_{s_1-6}$. 
We calculate
\[
P(u) - P(v) 
= \int_0^1 P'(v + \lm (u-v)) \, d\lm \, [u - v] =: \mL(u,v)[u - v]\ .
\]
The linear operator $\mL(u,v)$ has the same structure
as the operator $\mL_0$ in (2.12) of \cite{BFH}.
Since $u$ and $v$ both satisfy the Cauchy problem \eqref{i9}, 
we have $\mL(u,v)[u - v] = 0$ and $(u - v)(0) = 0$. 
Hence the well-posedness result in Lemma 6.7 of \cite{BFH}
implies $(u-v)(t) = 0$ for all $t \in [0,T]$.
This completes the proof of Theorem \ref{thm:control}.
\qed

\bigskip

\noindent
{\bf Proof of Theorem \ref{thm:byproduct}.}
We define 
\begin{align} \label{def Es bis}
E_s & := C([0,T], H^{s+6}_x) \cap C^1([0,T], H^{s+3}_x) 
\cap C^2([0,T], H^s_x),
\\
F_s & := \{ (g_1, g_2) : g_1 \in C([0,T], H^{s+6}_x) \cap C^1([0,T], H^s_x), 
g_2 \in H^{s+6}_x \}
\end{align}
equipped with norms 
\begin{align} \label{def norm Es bis}
\| u \|_{E_s} & := \| u \|_{T,s+6} + \| \pa_t u \|_{T,s+3} + \| \pa_{tt} u \|_{T,s}
\\
\| (g_1,g_2) \|_{F_s} & := \| g_1 \|_{T,s+6} + \| \pa_t g_1 \|_{T,s} + \| g_2 \|_{s+6},
\end{align}
and $\Phi(u) := (P(u), u(0))$, where $P$ is defined in \eqref{ge2}.
Given $g := (0,u_{in}) \in F_{s_0}$, 
the Cauchy problem \eqref{i11}
writes $\Phi(u) = g$.
We fix 
$V := \{ u \in E_3 : \| u \|_{E_3} \leq \d_0 \}$, 
where $\d_0$ is the same as in the proof of Theorem \ref{thm:control}; 
we fix $a_0, \mu, a_1, \a, \b,a_2$ like in \eqref{param.1},
where $\s$ is now the constant appearing in Lemma 6.7 of \cite{BFH},
$\t = \s + 9$ by Lemmas 2.1 and 6.7 of \cite{BFH} (combined with the definition of the spaces $E_s,F_s$),
$r \geq r_0 := a_2 + \t$ is the regularity of $\mN$,
and $\d_1$ is small enough to satisfy the assumption $\d(0) \leq \delta_*$ 
in Lemma 6.7 of \cite{BFH}.

Assumption \eqref{tame in NM} about the right inverse of the linearized operator 
is satisfied by Lemmas 6.7 and 2.1 of \cite{BFH}. 
We fix $s_0 := \a + 6$.
Then Theorem \ref{thm:1} applies, giving the existence part of Theorem \ref{thm:byproduct}. 
The uniqueness of the solution is proved exactly as in the proof of Theorem \ref{thm:control}.
This completes the proof of Theorem \ref{thm:byproduct}.
\qed

\begin{remark} \label{rem:better}
Although the linearized control problem \eqref{3009} admits a right inverse
with no loss of regularity in its argument (see \eqref{2809}, where $h,\ph$ have the same regularity $s$ as $g$),  
the application of H\oe{}rmander's implicit function theorem in Sobolev class 
gives a solution $f,u$ of the nonlinear control problem \eqref{i9}-\eqref{i10}
that is less regular, with arbitrarily small loss, than the data. 
This loss is due to the inclusion of the weak space $E_\a'$ into the spaces $E_a$ for all $a < \a$.
Thus, for initial and final states $u_{in}, u_{end} \in H^{s_1}$, 
the controllability theorem in \cite{BFH} (Theorem 1.1 of \cite{BFH}) 
gives the existence of a control 
$f$ and a solution $u$ of \eqref{i9}-\eqref{i10} 
of regularity 
\[
u,f \in C([0,T],H^s) \cap C^1([0,T], H^{s-3}) \cap C^2([0,T], H^{s - 6})
\quad \forall s < s_1,
\] 
with estimate 
\begin{multline*} 
\| u,f \|_{C([0,T],H^s_x)} + \| \pa_t u, \pa_t f \|_{C([0,T],H^{s-3}_x)} 
+ \| \pa_{tt} u, \pa_{tt} f \|_{C([0,T],H^{s-6}_x)} 
\\
\leq C_s (\| u_{in} \|_{s_1} + \| u_{end} \|_{s_1})
\quad \forall s < s_1,
\end{multline*}
for some constant $C_s > 0$, depending on $s,T,\om,\mN$,
and possibly diverging as $s \to s_1$.
The improvement of Theorem \ref{thm:control} with respect to the controllability theorem in \cite{BFH} is the achievement of the sharp, natural regularity $s_1$ of the problem,
without loss. 

Analogously, the improvement of Theorem \ref{thm:byproduct} 
with respect to the corresponding local existence and uniqueness theorem in \cite{BFH} 
for the Cauchy problem \eqref{i11} (Theorem 1.4 in \cite{BFH}) 
is the achievement of the sharp, natural regularity $s_0$ of the problem,
without loss 
(where ``sharp'' means that the solution has the same regularity as the datum).
\end{remark}

\begin{remark} \label{rem:3108}
The approach to control and Cauchy problems 
that we have used in the proof of Theorems \ref{thm:control} and \ref{thm:byproduct} 
also applies to other equations. 

In \cite{BHM} a similar result is proved for Hamiltonian, quasi-linear perturbations 
of the Schr\oe{}dinger equation on the torus in dimension one,
using Theorem \ref{thm:1}. 

Theorem \ref{thm:1} could also be used as an alternative approach, 
based on a different nonlinear scheme, 
to prove the controllability result for gravity capillary water waves in \cite{ABH}.

In the context of KAM for PDEs, 
Theorem \ref{thm:1} is used in \cite{BBHM-WW-h} 
to solve a quasi-periodic nonlinear PDE of the form 
$\om \cdot \pa_\ph u(\ph,x) = V(\ph, x + u(\ph,x))$ 
on the torus $(\ph,x) \in \T^{n+1}$, where $\om \in \R^n$ is a Diophantine vector.
This is the equation of the characteristic curves 
of a quasi-periodic transport equation.
\end{remark}

\begin{footnotesize}

\end{footnotesize}

\bigskip

\begin{flushright}

Pietro Baldi, Emanuele Haus 

\medskip

Dipartimento di Matematica e Applicazioni ``R. Caccioppoli''

Universit\`a di Napoli Federico II  

Via Cintia, 80126 Napoli, Italy

\medskip

\texttt{pietro.baldi@unina.it} 

\texttt{emanuele.haus@unina.it}
\end{flushright}

\end{document}